%% file: main.tex
\documentclass[ba,preprint]{imsart}

\RequirePackage{amsthm,amsmath,amsfonts,amssymb}
\RequirePackage[authoryear]{natbib}
\RequirePackage[colorlinks,citecolor=blue,urlcolor=blue,backref=page,backref=page]{hyperref}
\RequirePackage{graphicx}

\usepackage[utf8]{inputenc} 
\usepackage[T1]{fontenc}    
\usepackage{hyperref}       
\usepackage{url}            
\usepackage{booktabs}       
\usepackage{amsfonts}       
\usepackage{amsmath}
\usepackage{amsthm}
\usepackage[algoruled,vlined]{algorithm2e}
\usepackage{mathtools}
\usepackage{nicefrac}       
\usepackage{microtype}      
\usepackage{fancyhdr}       
\usepackage{graphicx}       
\usepackage[export]{adjustbox}
\usepackage{aligned-overset}
\graphicspath{{media/}}     
\usepackage{enumitem}

\usepackage{xspace}
\usepackage{xcolor}

\allowdisplaybreaks    

\let\originalleft\left
\let\originalright\right
\renewcommand{\left}{\mathopen{}\mathclose\bgroup\originalleft}
\renewcommand{\right}{\aftergroup\egroup\originalright}

\newcommand{\matern}{Mat\'ern\xspace}

\newcommand{\bm}{\mathbf{m}}

\newcommand{\bx}{\mathbf{x}}
\newcommand{\by}{\mathbf{y}}
\newcommand{\btheta}{{\boldsymbol\vartheta}}
\newcommand{\bxi}{{\boldsymbol\xi}}

\newcommand{\calG}{\mathcal{G}}
\newcommand{\calN}{\mathcal{N}}

\newcommand{\bbE}{\mathbb{E}}
\newcommand{\bbR}{\mathbb{R}}

\pubyear{2024}
\volume{TBA}
\issue{TBA}
\firstpage{1}
\lastpage{1}

\startlocaldefs
\theoremstyle{plain}
\newtheorem{theorem}{Theorem}[section]
\newtheorem{lemma}[theorem]{Lemma}
\newtheorem{corollary}[theorem]{Corollary}
\newtheorem{proposition}[theorem]{Proposition}

\theoremstyle{definition}

\newtheorem{assumption}{Assumption}
\newtheorem{remark}{Remark}

\def\disclosurename{Disclosure of tool use}
\newenvironment{disclosure}[1][\disclosurename]{\begin{acks}[\disclosurename]}{\end{acks}}
\endlocaldefs

\begin{document}

\begin{frontmatter}
\title{Mitigating surrogate model error in inverse problems via generalised Bayesian inference}
\runtitle{A sample running head title}

\begin{aug}
\author[A]{\fnms{Hans}~\snm{Reimann}\ead[label=e1]{reimann@math.uni-heidelberg.de}},
\author[B]{\fnms{Simon}~\snm{Urbainczyk}\ead[label=e2]{}},
\author[C]{\fnms{Jonas}~\snm{Latz}\ead[label=e3]{jonas.latz@manchester.ac.uk}},
\and
\author[D]{\fnms{Aretha L}~\snm{Teckentrup}\ead[label=e4]{a.teckentrup@ed.ac.uk}}


\address[A]{Institute for Mathematics,
Heidelberg University\printead[presep={,\ }]{e1}}
\address[B]{Maxwell Institute and School of Mathematical and Computer Sciences,
Heriot-Watt University\printead[presep={ \ }]{e2}}
\address[D]{Department of Mathematics,
University of Manchester\printead[presep={,\ }]{e3}}
\address[C]{Maxwell Institute and School of Mathematics,
University of Edinburgh\printead[presep={,\ }]{e4}}
\runauthor{H. Reimann et al.}
\end{aug}

\begin{abstract}
Employing surrogate forward operators in place of the true mathematical model in inverse problems significantly improves computational efficiency, yet induces approximation error causing potentially unreliable inference. We explore generalised Bayesian inference with weighted cross-entropy to adjust uncertainty quantification in inverse problems, with a focus on Gaussian process surrogate models. The proposed modular approach provides theoretical guarantees in global bias robustness and establishes posterior accuracy via decomposition of surrogate approximation and noisy observation errors. For the considered class of Gaussian process surrogate models, we obtain  efficient sampling methods for the generalised posteriors via closed form gradients and a novel tuning strategy for the learning rate. The proposed approaches are evaluated in simulation studies on an inverse problem in partial differential equations, investigating uncertainty quantification and degree of misspecification. Our results contribute to theoretical understanding of generalised posteriors in Bayesian inverse problems as well as implications for algorithmic practice. 
\end{abstract}



\end{frontmatter}

\section{Introduction}

Bayesian inverse problems (BIPs) concern Bayesian inference in the context of mathematical models that may describe e.g. a physical system. These mathematical models usually enter the inferential procedure as sophisticated observation (or forward) operators that are part of the likelihood, resulting in likelihoods that are very costly to evaluate numerically. A typical example of an underlying model are the partial differential equations (PDEs) found in numerous applications in the natural sciences, engineering, and medicine.

For uncertainty quantification (UQ) via Bayesian inference, sampling procedures for exploring the posterior distribution typically need to evaluate the forward operator, and often its gradient, many times. In these situations, with severe constraints on numerical complexity, employing cheap, pre-trained surrogate models (or {\em emulators, reduced order models}) for the forward operator provides a valuable reduction in cost. Gaussian process (GP) surrogates \citep{sacks1989design,kennedy2001bayesian,williams2006gaussian,lykkegaard2023multilevel,helin2024introduction} are an efficient and popular choice of model to approximate the action of the forward operator. Alternatively, neural network surrogates have also surged in popularity in recent years \citep{deveney2023deep,jimenez2026neural}.

While exploiting a surrogate model yields a significant reduction in computational cost, it also introduces a source of error via its approximate nature. In practice, where the accuracy of the surrogate model may be limited, this error can be significant and lead to biased and overconfident inference. To mitigate this issue, recent works including \citep{conrad2017statistical,calvetti2018iterative,fox2020randomized,cleary2021calibrate,lykkegaard2023multilevel,bai2024gaussian} have advocated for incorporating a model for the surrogate error explicitly in the approximate likelihood. The error model is chosen as a Gaussian distribution, either coming directly from using the posterior variance of a Gaussian process surrogate \citep{calvetti2018iterative,bai2024gaussian} or being estimated for a surrogate based on a coarse numerical solver \citep{conrad2017statistical,fox2020randomized}. The concept is closely related to the inclusion of model discrepancy in Bayesian inference of \citep{kennedy2001bayesian,brynjarsdottir2014learning}, where a similar error model is used to mitigate the error between a mathematical model and a real world process.

The assumption that the surrogate error follows a Gaussian distribution may be unjustified or difficult to verify in practice, and thus, an approximate likelihood in which a GP surrogate replaces a complex forward operator may be misspecified even after accounting for the surrogate error through the posterior variance.
We therefore investigate generalised Bayesian inference (GBI) to address misspecification of an assumed Gaussian observation likelihood in BIPs. The corresponding generalised posterior provides surrogate error robustness by controlling highly influential instances of departure from the assumed Gaussian observation model, mitigating the unknown epistemic error of the surrogate models. We refer to \citet{bissiri2016general}, \citet{jewson2018principles}, and \citet{knoblauch2019generalized} for pioneering works investigating generalised Bayesian inference and \citet{boustati2020generalised}, \citet{duran2024outlier}, and \citet{reimann2026kalman} for utilising generalized posteriors in the context of BIPs. GBI and generalised posteriors are further closely related to Gibbs posteriors in PAC-Bayesian learning (see \citealp{mcallester1998some}), ABC approaches (see \citealp{del2012adaptive,lee2012choice}), and likelihood-free inference (see \citealp{sisson2011likelihood}) in their respective contexts.

\paragraph{Setup and motivation.}
Although we motivate investigation via emulators in inverse problems involving solution operators of PDEs, our approach applies more general to observation likelihoods at risk of misspecification. We consider the inverse problem of recovering a set of true parameters $\btheta^\dagger\in T \subset \bbR^{d_{\btheta}}$ given implicitly through
\begin{equation}
     \by=\calG(\btheta^\dagger) + {\boldsymbol\eta}^\dagger
    \label{eq:obs_ip}
\end{equation}
where $\by\in\bbR^{d_\by}$ are noisy observations produced by a parameter-to-observation map $\calG: T \to \bbR^{d_\by}$ and ${\boldsymbol\eta}^\dagger$ is a realisation of the observational noise. We generally assume Gaussian observational noise ${\boldsymbol\eta} \sim \mathcal{N}(0, \Gamma_{\boldsymbol\eta})$, often for the sake of numerical tractability and limited knowledge. Note, misspecification regarding $\boldsymbol{\eta}$ e.g. in tail decay strongly impacts inference. While the proposed methods do also provide robustness to misspecification in the observational noise, this is not the focus of the current work.

The Bayesian approach to inverse problems additionally assumes an a priori probability density  $p_0(\btheta)$ on $\btheta\in T$. Via Bayes' rule and after after observing $\by$, we obtain the posterior density 
\begin{equation}
    \label{eq:posterior}
    p(\btheta \, | \, \by)\propto \ell(\by \, | \, \btheta) \, p_0(\btheta)
\end{equation}
with $\ell(\by \, | \, \btheta)$ denoting the likelihood of observing $\by$ given parameter $\btheta$. Assuming Gaussian observational noise ${\boldsymbol\eta} \sim \mathcal{N}(0, \Gamma_{\boldsymbol\eta})$, we have
\begin{equation}
    \label{eq:obs_likelihood}
    \ell(\by \, | \, \btheta) \propto \exp\left(-\frac{1}{2} \left\|\calG(\btheta) - \by\right\|^2_{\Gamma_{\boldsymbol\eta}}\right)\quad\text{with}\quad\|\bx\|^2_{\Gamma_{\boldsymbol\eta}}=\left\langle\bx,\,\Gamma_{\boldsymbol\eta}^{-1}\bx\right\rangle\,.
\end{equation}

We aim to avoid the (potentially costly) evaluation of the observation operator $\calG$ by employing a surrogate, in particular a GP based surrogate, and hence altering the likelihood term. Let $\widetilde{\calG}:T\to\mathbb{R}^{d_\by}$ be a pre-trained surrogate model emulating the action of the parameter-to-observation map $\calG$. Replacing $\calG$ by $\widetilde{\calG}$ in \eqref{eq:obs_ip} introduces an additional, epistemic approximation error. Thus, we obtain
\begin{equation*}
        \by= \widetilde{\calG}(\btheta^\dagger) + {\boldsymbol\zeta}_\btheta^\dagger + {\boldsymbol\eta}^\dagger\,.
\end{equation*}
The approximation term $\boldsymbol{\zeta}_\btheta^\dagger\coloneqq\calG(\btheta^\dagger) - \widetilde{\calG}(\btheta^\dagger)$ introduces a bias depending on the true parameter $\btheta^\dagger$. While not inherently random, we model it as a random perturbation with support in $\mathbb{R}^{d_\by}$ w.r.t. the prior. The observational error then combines the aleatoric stochastic error ${\boldsymbol\eta}^\dagger$ and an epistemic approximation error ${\boldsymbol\zeta}_\btheta^\dagger$, and is generally unknown in structure. The model in \eqref{eq:posterior} is no longer well-specified when employing a surrogate model and approximate Bayesian inference methods can improve accuracy (see the discussion in \citealp{knoblauch2022optimization}). Yet, the Gaussian observational error assumed in \eqref{eq:obs_likelihood} may still provide a best guess and maintains tractability of the inverse problem. For negligible instances of epistemic approximation error ${\boldsymbol\zeta}_\btheta^\dagger$ regular Bayesian methods remain accurate. In the contrary instance of malicious error ${\boldsymbol\zeta}_\btheta^\dagger$, they produce overconfident inference. Generalized posteriors can be used to target influential instances of approximation error ${\boldsymbol\zeta}_\btheta^\dagger$ and maintain computational feasibility at the cost of optimality of the Bayesian posterior distribution in well-specified cases.

\paragraph{Relevance, contribution and structure.}
We expand recent advances in generalized Bayesian inference to common practice in Bayesian inverse problems with non-linear forward operators and surrogates models. Although our approach is motivated from a perspective of robustness to model-observation mismatch, it applies to the wider context of robust Bayesian inference under observation likelihood misspecification. Specifically, our contributions are
\vspace{-0.4cm}
\begin{itemize}
    \item a modular design framework combining generalized posteriors via an adjusted weighted cross-entropy and Gaussian process emulators of the forward operator,
    \vspace{-0.20cm}
    \item deriving conditions for posterior global bias robustness for non-linear observation operators and establishing a robustness-consistency trade-off in posterior Hellinger distance for insights into accuracy, and
    \vspace{-0.20cm}
    \item addressing practical aspects via posterior gradients and heuristic tuning and a simulation study on an inverse problem of a PDE in one dimension for illustration and insights.
    \vspace{-0.3cm}
\end{itemize}
We first introduce Gaussian process emulators in Bayesian inverse problems and generalised Bayesian inference via weighted cross-entropy to create a modular framework of specifications in Section~\ref{sec:background}. We proceed with a theoretical analysis in global bias robustness, and accuracy and consistency in Section~\ref{sec:the_prop}. Next, we consider practical implementation via closed form expressions of the generalized posterior score function within the modular framework and tuning of a novel degree of freedom in Section~\ref{sec:sampling_methodology}. We illustrate the proposed methods in a simulation experiment presented in Section~\ref{sec:numerical_experiment} before discussion and conclusions in Section~\ref{sec:disc_conclusion}.

\section[Generalised BIPs with GP surrogates]{Generalised Bayesian inversion with Gaussian process surrogates}

\label{sec:background}
We combine the ideas of GP surrogates and GBI to obtain an alternative to classical Bayesian inverse problems that remains computationally efficient yet controls the influence of the induced surrogate approximation error via bias robustness. The result is a modular framework in choice for GP surrogate and parametrisation in GBI with weighted cross-entropy.

\subsection{Gaussian processes surrogates}

For evaluation of the posterior in \eqref{eq:posterior} via sampling, one typically requires a large number of costly forward solves. We utilise the class of GP surrogates as an affordable approximation to the exact likelihood. 

We consider GPs $\mathcal{GP}(m, K)$ with input dimension $d_{\btheta}$, output dimension $d_\by$, mean function $m$, and covariance kernel $K(\btheta, \btheta')=k(\btheta, \btheta')\,\mathrm{Id}_{d_\by}\in\bbR^{d_\by \times d_\by}$ for positive definite (scalar) kernel $k$ and $\btheta,\btheta'\in T$. A common choice of covariance kernel is the \matern kernel (see \citealp{williams2006gaussian}) recovering a squared exponential kernel in the smoothness limit. While we model the individual GP output dimensions independent for simplicity, additional knowledge such as spatial correlation can be encoded in the covariance kernel (see e.g. \citealp{bai2024gaussian,raissi2017machine,chkrebtii2016,cockayne2017,pfortner2022physics}).

For approximating the action of $\calG$, we assume access to (noisy or noise-free) evaluations of $\calG$ at a set of $n\in\mathbb{N}$ distinct training locations in $T$ and place a GP prior $\mathcal{GP}(0, K)$ on $\calG$. We update the prior via conditioning on the training evaluations of $\calG$ to obtain the GP posterior $\widetilde{\calG}_n$ with mean $\bm_n^{\calG}(\btheta)$ and covariance $K_n^{\calG}(\btheta, \btheta')$ via regression \citep{williams2006gaussian}. We utilise the GP posterior $\widetilde{\calG}_n$ for predicting $\calG$ in the likelihood term \eqref{eq:obs_likelihood}.

Agnostic to the randomness of $\widetilde{\calG}_n$, the action of $\calG$ may be approximated by the posterior mean $\bm_n^{\calG}$ \citep{stuart2018posterior,bai2024gaussian}. In this case we approximate $\ell$ by
\begin{equation}
    \label{eq:mean_surrogate}
    \widetilde{\ell}_\mathrm{mean}(\by \,|\, \btheta)
    \propto \exp\left(-\frac{1}{2} \left\|\bm_n^{\calG}(\btheta) - \by\right\|^2_{\Gamma_{\boldsymbol\eta}}\right)\,.
\end{equation}
This approximation has been observed to result in overconfident posteriors since the covariance of the BIP posterior fails to capture the additional uncertainty induced by the epistemic approximation error of the GP mean \citep{bai2024gaussian}. Note, this essentially reduces to emulation via kernel interpolation \citep{bishop2006pattern}.

As a remedy, the GP marginal approximation was proposed. Via employing information on covariance of the GP emulator of $\calG$, the error induced by the surrogate can be partially mitigated \citep{bai2024gaussian,calvetti2018iterative,stuart2018posterior}. Following \citet[Appendix A]{bai2024gaussian}, we obtain the approximate likelihood
\begin{align}
    \label{eq:marginal_surrogate}
    \widetilde{\ell}_\mathrm{GP}(\by \,|\, \btheta)
    &\propto \bbE_{\widetilde{\calG}_n(\btheta)}\left[\exp\left(-\frac{1}{2} \left\|\widetilde{\calG}_n(\btheta) - \by\right\|^2_{\Gamma_{\boldsymbol\eta}}\right)\right] \nonumber \\
    &\propto \frac{1}{\sqrt{\det(K_n^{\calG}(\btheta,\btheta) + \Gamma_{\boldsymbol\eta})}} \, \exp\left(-\frac{1}{2} \left\|\bm_n^{\calG}(\btheta) - \by\right\|^2_{\left(K_n^{\calG}(\btheta,\btheta) + \Gamma_{\boldsymbol\eta}\right)}\right)\,.
\end{align}
The additional uncertainty considered by the marginal likelihood can improve UQ, but typically fails to fully capture the structure and effect of the surrogate approximation error. The performance of either posterior approximation in practice strongly depends on the number of training points $n$ of the GP surrogate, as well as the choice of prior distribution \citep{bai2024gaussian}.
This motivates approaches to BIPs robust to model-observation mismatch.

\subsection{Generalised Bayesian inference}\label{sec:gen_Bayes}

Generalised Bayesian inference is concerned with approximate posterior distributions aiming to improve inference when assumptions are not met. Misspecification of observation models provides an influential instance of such (see e.g. \citealp{knoblauch2019generalized}). 

Regular Bayesian inference is information optimal \citep{zellner1988optimal}, with information in the sense of \citet{kullback1951information}. We observe that
\begin{align}
    \label{eq:KL_Bayes}
    p(\btheta \,|\, \by)&=\frac{\ell(\by\,|\, \btheta)p_0(\btheta)}{p(\by)}=\exp\left(-\widehat{\mathrm{KL}}\left[p(\by) \,\|\, \ell(\by \,|\, \btheta)\right]\right)\, p_0(\btheta) \nonumber \\
    &\propto \ell(\by\,|\, \btheta)\, p_0(\btheta)=\exp\left(-\widehat{\mathrm{CE}}\left[p(\by) \,\|\, \ell(\by \,|\, \btheta)\right]\right)\, p_0(\btheta)
\end{align}
with $\widehat{\mathrm{KL}}\left[p(\by) \,\|\, \ell(\by \,|\, \btheta)\right] = \log\left(\frac{p(\by)}{\ell(\by|\btheta)}\right)$ the one-sample estimator of the Kullback-Leibler (KL) divergence and $\widehat{\mathrm{CE}}\left[p(\by) \,\|\, \ell(\by \,|\, \btheta)\right] = -\log\ell(\by\,|\, \btheta)$ the one-sample cross-entropy (CE) estimator. This assumes $\by$ is sampled from the evidence $p(\by)=\int_T\ell(\by\,|\, \btheta)p_0(\btheta)d\btheta$.

When components such as the observation likelihood are misspecified and the evidence does not represent the true data generating process in relevant properties such as tail decay, theoretical guarantees including optimal information processing may no longer hold. Regular posteriors are known to be sensitive to outliers \citep{knoblauch2019generalized,pacchiardi2024generalized} and even more so for small numbers of observations. This mainly follows from properties of the KL divergence.

Generalized Bayesian inference replaces the KL divergence in \eqref{eq:KL_Bayes} by an alternative discrepancy measure $\mathrm{D}$ of distributions on $\mathbb{R}^{d_\by}$. The corresponding generalized posterior is then
\begin{align}\label{eq:generalised_posterior}
    p_{\beta\mathrm{D}}(\btheta \,|\, \by)
    \propto \exp\left(-\beta\widehat{D}\left[p(\by) \,\|\, \ell(\by \,|\, \btheta)\right]\right) \, p_0(\btheta)\,,
\end{align}
with $\widehat{\mathrm{D}}$ a corresponding one sample estimator and learning rate (or inverse temperature) $\beta>0$. The case $\mathrm{D}=\mathrm{KL}$ and $\beta=1$ recovers \eqref{eq:KL_Bayes}.
The choice of discrepancy measure $\mathrm{D}$ may target violation of assumptions while subject to restrictions in applicability. Early work targeting observation likelihood misspecification in different contexts via robust posteriors \citep{ghosh2016robust,knoblauch2018spatio,boustati2020generalised} considered $\beta$- and power divergences with more recent work in \citet{altamirano2023robustBOCPD,altamirano2023robustGPR} and \citet{reimann2026kalman} employing diffusion score matching. 

We investigate an adjusted weighted cross-entropy (WCE) utilised in \citet{duran2024outlier} for linear BIPs in Kalman filtering. Let
\begin{equation}\label{eq:wce}
    \mathrm{WCE}\left[p(\by) \,\|\, \ell(\by \,|\, \btheta)\right]
    \coloneqq\bbE_{p}\left[- w^2(\by, \btheta) \log \ell_+(\by \,|\, \btheta)\right]
\end{equation}
where $w:\mathbb{R}^{d_\by}\times T\to \bbR$ is a suitable weight function and $\ell_+$ is the observation likelihood $\ell$ normalised such that
\begin{equation*}
    -\log \ell(\by \,|\, \btheta)=-\log\ell_+(\by \,|\, \btheta)+C(\btheta)
    \quad\text{with}\quad
    C(\btheta)\coloneqq\inf_{\by'\in\bbR^{d_\by}}\left[-\log\ell(\by'\,|\,\btheta)\right]
\end{equation*}
and $-\log\ell_+(\by \,|\, \btheta) \geq 0$ with equality at the mode. Our choice is motivated by the properties that WCE exhibits an immediate one-sample estimator and guarantees posterior robustness for suitable choice of weight $w$ (see Assumption~\ref{ass: kernel} in Section~\ref{sec:the_prop}) via likelihood rescaling based on model-observation mismatch. 

\begin{remark}\label{rem:wce}
    We replace the log likelihood by $\ell_+$ to avoid scaling an observation independent component $C(\btheta)$ by $w^2(\by, \btheta)$. Furthermore, we also drop $C(\btheta)$ such that the generalised posterior recovers the prior $p_0$ in the limit $w\to0$ instead of the tilted term $\exp\left(-\beta C(\btheta)\right) \, p_0(\btheta)$. For weights $w$ and $C$ both independent of $\btheta$, the adjusted WCE recovers the generalised posterior in \citet{duran2024outlier}. Assuming Gaussian observation likelihoods as in \eqref{eq:obs_likelihood} yields $-\log\ell_+(\by\,|\,\btheta)=\frac{1}{2}\|\calG(\btheta)-\by\|^2_{\Gamma_{\boldsymbol\eta}}$, while the GP marginal approximation in \eqref{eq:marginal_surrogate} loses the parameter dependent normalizing constant $\det(K_n^{\calG}(\btheta,\btheta) + \Gamma_{\boldsymbol\eta})^{-\frac{1}{2}}$. 
\end{remark}

When choosing the weight $w$ as a function of $(\by,\btheta)\in\bbR^{d_\by}\times T$, we refer to it as \emph{dynamic weighting} with respect to the parameter of interest. Alternatively, when choosing $w$ as a function of $\by\in\bbR^{d_\by}$ only or replacing $\btheta\in T$ by a fixed value, we refer to it as \emph{static weighting}. A default choice then is the prior mean $\bar\btheta = \bbE_{ p_0}[\btheta]$, so $w^2(\by,\bar\btheta)$.

Substituting the WCE estimator into the generalised posterior \eqref{eq:generalised_posterior}, the learning rate $\beta$ directly scales the range of the weight $w$. Moreover, $\beta w^2$ essentially combines into a single parameter and observation dependent learning rate. For choice of $\beta=1$, WCE, and $w$ with range $(0,1]$, the corresponding generalized posterior bridges between the prior and the regular Bayesian posterior. While this is desirable to counteract overconfidence, we also investigate a learning rate $\beta=2$ to offset excessive uncertainty inflation.

As in \citet{altamirano2023robustBOCPD,altamirano2023robustGPR} and \citet{duran2024outlier}, we choose $w$ as an adjusted inverse multi-quadratic (IMQ) function
\begin{equation}\label{eq:w_squared}
    w^2(\by, \btheta)
    = \Bigg(1 + \frac{\|\widetilde{\calG}(\btheta) - \by\|^2_{\Gamma(\btheta)}}{q^2}\Bigg)^{-1}
\end{equation}
with parameter $0<q^2<\infty$ and $\widetilde{\calG}(\btheta)$ and $\Gamma(\btheta)$ to match the likelihood $\ell_+$ considered by CE term, e.g. $\calG(\btheta)$ and $\Gamma_{\boldsymbol\eta}$ in \eqref{eq:obs_likelihood}. The IMQ kernel satisfies assumptions for robustness in Section~\ref{sec:the_prop}.

\subsection{Pick and choose}
We introduced two concepts that make changes to the assimilation of an observation in Bayesian inverse inference for reasons that can be seen as perpendicular. Surrogate forward models serve to reduce computational cost by replacing the action of the forward operator with a cheaper approximation. Generalized posteriors can counteract deviation from assumptions, e.g. the observational error deviating from a Gaussian assumption. This work proposes to combine both concepts. The epistemic approximation error induced by emulation contaminates the assumed observation likelihood. The choices of GP surrogates of the forward operator and WCE-based GBI maintain closed form posterior score functions. An overview of the resulting modular framework of WCE posteriors with choice in GP surrogate, learning rate, and weight parametrisation is provided in Table~\ref{tab:props}.

\begin{table}[tb]
    \centering
    \begin{tabular}{lll}
        \toprule
        \textbf{GP surrogate likelihood} & \textbf{Learning rate} & \textbf{Weight parametrisation} \\
        \midrule
        mean, $\widetilde{\ell}_\mathrm{mean}$ & $\beta=1$ & static, $w^2(\by,\bar\btheta)$ \\
        marginal, $\widetilde{\ell}_\mathrm{GP}$ & $\beta=2$ & dynamic, $w^2(\by,\btheta)$ \\
        \bottomrule
        \addlinespace[\belowrulesep]
    \end{tabular}
    \caption{Modular framework of choices for WCE posteriors with GP surrogate forward operators. Each column provides a component with $8$ possible variants in total.}
    \label{tab:props}
\end{table}

We denote approximate likelihoods by $\widetilde{\ell}$, and WCE posteriors with dynamic weights by $\widetilde{p}_\beta$ and with static weights by $\bar{p}_\beta$, so
\begin{align*}
    \widetilde{p}_{\beta}(\btheta \,|\, \by)
    &\propto \exp\left(\beta w^2(\by, \btheta) \log \widetilde{\ell}_+(\by \,|\, \btheta)\right) \, p_0(\btheta)\quad\text{and}\\
    \bar{p}_{\beta}(\btheta \,|\, \by)
    &\propto \exp\left(\beta w^2(\by, \bar\btheta) \log \widetilde{\ell}_+(\by \,|\, \btheta)\right) \, p_0(\btheta)\,.
\end{align*}

\section{Theoretical properties}\label{sec:the_prop}

We investigate two theoretical aspects of the proposed approaches. First, we establish global bias robustness via a uniform bound on the posterior influence function. Second, we explore posterior accuracy and consistency via a trade-off between epistemic and aleatoric error. We aim to keep statements general and consider arbitrary approximations $\widetilde{\calG}$ of the true parameter-to-observation map $\calG$. We focus on $ \widetilde{p}_{\beta}$ with $\bar{p}_{\beta}$ considered as a special case. The global bias robustness directly applies for the mean and marginal GP surrogates, while the analysis on accuracy considers surrogates such as the mean GP surrogate first and later generalizes results to the marginal GP surrogate. 

\subsection{Global bias robustness}\label{sec:gbr}

We approach to posterior robustness as in \citet{altamirano2023robustBOCPD,altamirano2023robustGPR} and make explicit use of \citet[Proposition B.1]{altamirano2023robustBOCPD}. The weight $w$ introduces an additional degree of freedom and ensures robustness via its properties. Global bias robustness of posterior distributions closely relates to robustness in the sense of \citet{huber2004robust} and gross-error sensitivity in a sense of \citet{hampel2011robust}. Yet, as we consider posterior distributions as estimators of distributions, the classical concept of the influence function for measuring effect of contamination on point estimators is adjusted to the posterior influence function (PIF).

We model the true data generating distribution (DGD) $\pi^\dagger$ as an $\varepsilon$-contamination of the assumed data distribution $\widetilde{\pi}$ associated with the assumed observational noise and approximate forward operator combined in $\widetilde{\ell}$. The true DGD enters the one sample estimator of the WCE in \eqref{eq:wce} only via an observation $\by\in \bbR^{d_\by}$. Write equivalently
\begin{equation*}
     p_{\beta\mathrm{WCE}}(\btheta \,|\, \pi^\dagger)
    \propto \exp\left(-\beta\,\widehat{\mathrm{WCE}}\left[\pi^\dagger \,\|\, \ell(\,\cdot \,|\, \btheta)\right]\right) \, p_0(\btheta)\,.
\end{equation*}
For a Dirac measure $\pi^\dagger=\delta_\by$ this recovers exactly the one-sample WCE estimator for a single observation $\by\in \bbR^{d_\by}$ of the true DGD. We consider the $\varepsilon$-contamination
\begin{equation*}
    \pi^\dagger=\pi^\dagger_{\varepsilon,\by}=(1-\varepsilon)\widetilde{\pi}+\varepsilon\delta_{\by}\quad\text{with}\quad\by\in\bbR^{d_\by}\ \text{and}\ \varepsilon\in[0,1]\,. 
\end{equation*}
The contaminated measure $\pi^\dagger$ is the single instance departure of the true DGD from an assumed data distribution $\widetilde{\pi}$. The posterior influence function is then
\begin{equation*}
    \mathrm{PIF}(\by,\btheta,\widetilde{\pi})=\left.\frac{\mathrm{d}}{\mathrm{d}\varepsilon} p_{\beta\mathrm{WCE}}(\btheta|\pi^\dagger_{\varepsilon,\by})\right\vert_{\varepsilon=0} \,,
\end{equation*}
a point-wise expansion of the regular influence function over the parameter space $T$. We consider the posterior globally bias robust if the PIF is uniformly bounded:
\begin{equation}\label{eqn: global bias rob}
    \sup_{\btheta\in T,\ \by\in\mathbb{R}^{d_\by}}\mathrm{PIF}(\by,\btheta,\widetilde{\pi})<\infty \,.
\end{equation}
Regular posteriors with assumed Gaussian observation likelihoods are not robust in general (see \citealp[Theorem 3.2]{duran2024outlier} for the linear BIP case). We additionally refer to the works in \citet{ghosh2016robust} and \citet{matsubara2022robust} for investigations on robust posteriors and the posterior influence function. 

We generalize the result in \citet{duran2024outlier}, which assumed unadjusted WCE, linear observation operators, and Gaussian priors, to the given context. The stated result covers approximation error as well as observational noise misspecification. 

\begin{assumption}\label{ass: kernel}
For an assumed conditionally Gaussian likelihood 
\begin{equation*}
    \widetilde{\ell}_+(\by \, | \, \btheta)=\exp\left(-\frac{1}{2} \big\|\widetilde{\calG}(\btheta) - \by\big\|^2_{\Gamma(\btheta)}\right)\,,
\end{equation*}
let the weight function $w:\mathbb{R}^{d_\by}\times T\to (0,1]$ be such that $w^2(\by, \btheta)\big\|\widetilde{\calG}(\btheta)-\by\big\|^2_{\Gamma(\btheta)}\leq C(w)$ for all $(\by,\btheta)\in\mathbb{R}^{d_\by}\times T$ with $C(w)<\infty$.
\end{assumption}

The IMQ weight in \eqref{eq:w_squared} satisfies Assumption~\ref{ass: kernel} with $C(w)=q^2$ and $0<q^2<\infty$.

\begin{assumption}\label{ass: prior}
Let the prior distribution $p_0$ be such that $\sup_{\btheta\in T}p_0(\btheta)<\infty$.
\end{assumption}

\begin{theorem}[Global bias robustness]\label{thm:WCEgbr}
    Given the Bayesian inverse problem with assumed conditionally Gaussian likelihood 
    \begin{equation*}
        \widetilde{\ell}_+(\by \, | \, \btheta)= \exp\left(-\frac{1}{2} \big\|\widetilde{\calG}(\btheta) - \by\big\|^2_{\Gamma(\btheta)}\right)\,.
    \end{equation*} 
    Assume the prior $p_0$ satisfies Assumption~\ref{ass: prior} and the weight $w$ satisfies Assumption~\ref{ass: kernel}. The dynamic weighted cross-entropy posterior $p_{\beta\mathrm{WCE}}(\btheta \,|\, \widetilde{\pi})$ with $\beta\in\{1,2\}$ is globally bias robust in the sense of \eqref{eqn: global bias rob} for any probability measure $\widetilde{\pi}$.
\end{theorem}

\begin{proof}
    We employ \citet[Proposition B.1]{altamirano2023robustBOCPD} and show
    \vspace{-0.4cm}
    \begin{enumerate}
        \item $\sup_{\by\in\mathbb{R}^{d_\by}}\ \left|\beta\,\widehat{\mathrm{WCE}}\left[\delta_\by \,\|\, \widetilde{\ell}(\,\cdot \,|\, \btheta)\right]\right|<\gamma(\btheta)$ for $\gamma:T\to[0,\infty)$,
        \vspace{-0.2cm}
        \item $\sup_{\btheta\in T}p_0(\btheta)\gamma(\btheta)<\infty$, and $\int_Tp_0(\btheta)\gamma(\btheta)\mathrm{d}\btheta<\infty$.
        \vspace{-0.3cm}
    \end{enumerate}

    Regarding condition 1, we observe for any $\by\in\mathbb{R}^{d_\by}$ and via Assumption~\ref{ass: kernel} that
    \begin{align*}
    \left|\beta\,\widehat{\mathrm{WCE}}\left[\delta_\by \,\|\, \widetilde{\ell}(\,\cdot \,|\, \btheta)\right]\right|&= \left|\beta w^2(\by, \btheta) \log \ell_+(\by \,|\, \btheta)\right|\\
    &\leq \frac{2w^2(\by, \btheta)}{2}\big\|\widetilde{\calG}(\btheta) - \by\big\|^2_{\Gamma(\btheta)}\leq C(w)\,.
    \end{align*}
    
   Via Assumption~\ref{ass: prior} condition 2 holds since 
   \begin{equation*}
       \sup_{\btheta\in T}C(w)p_0(\btheta)<\infty\quad\text{and}\quad\int_TC(w)p_0(\btheta)d\btheta=C(w)\,.
   \end{equation*}
\end{proof}

We require additional assumptions for a static weight. 

\begin{corollary}\label{cor:statgbr}
Given the assumptions in Theorem~\ref{thm:WCEgbr},
\vspace{-0.3cm}
\begin{itemize}
    \item $w^2(\by,\btheta_0)\left\| \by\right\|^2_{\Gamma(\btheta)}\leq C_s(w)$ and $w^2(\by,\btheta_0)\big\langle \by,\Gamma^{-1}(\btheta)\widetilde{\calG}(\btheta)\big\rangle\leq\gamma_s(\btheta)$\\ for all $\by\in\bbR^{d_\by}$ with $C_s(w)<\infty$ and $\gamma_s:T\to[0,\infty)$, and
    \vspace{-0.1cm}
    \item $\sup_{\btheta\in T}p_0(\btheta)\big\|\widetilde{\calG}(\btheta)\big\|^2_{\Gamma(\btheta)},\ \sup_{\btheta\in T}p_0(\btheta)\gamma_s(\btheta)<\infty$ and\\ $\int_T p_0(\btheta)\big\|\widetilde{\calG}(\btheta)\big\|^2_{\Gamma(\btheta)}\mathrm{d}\btheta,\ \int_T p_0(\btheta)\gamma_s(\btheta)\mathrm{d}\btheta<\infty$.
    \vspace{-0.3cm}
\end{itemize}
The \emph{static weighted} cross-entropy posterior 
with $w^2(\by,\btheta_0)$ for an arbitrary, fixed $\btheta_0\in T$ is globally bias robust in the sense of \eqref{eqn: global bias rob} for any probability measure $\widetilde{\pi}$.
\end{corollary}

\begin{proof}
   The weight $w$ reduces to a function in $\by\in\mathbb{R}^{d_\by}$. We adjust the proof of Theorem~\ref{thm:WCEgbr} by splitting the quadratic form via
\begin{equation*}
    \big\|\widetilde{\calG}(\btheta) - \by\big\|^2_{\Gamma(\btheta)}=\big\|\widetilde{\calG}(\btheta)\big\|^2_{\Gamma(\btheta)}-2 \left\langle \by,\Gamma^{-1}(\btheta)\widetilde{\calG}(\btheta)\right\rangle+\left\| \by\right\|^2_{\Gamma(\btheta)}\,. 
\end{equation*}
Then for all $y\in\mathbb{R}^{d_\by}$
\begin{equation*}
    \left|\beta\,\widehat{\mathrm{WCE}}\left[\delta_\by \,\|\,\widetilde{\ell}(\,\cdot \,|\, \btheta)\right]\right|\leq\big\|\widetilde{\calG}(\btheta)\big\|^2_{\Gamma(\btheta)} + 2\gamma_s(\btheta) + C_s(w) \eqqcolon\gamma(\btheta)
\end{equation*}
via the first additional assumption. The second additional assumption ensure finiteness in
\begin{equation*}
    \sup_{\btheta\in T}p_0(\btheta)\gamma(\btheta)<\infty\quad\text{and}\quad\int_Tp_0(\btheta)\gamma(\btheta)\mathrm{d}\btheta<\infty\,.
\end{equation*}
\end{proof}

The additional assumptions place constraints in balancing prior and forward operator. When assuming $\sup_{\btheta\in T}\big\|\widetilde{\calG}(\btheta)\big\|^2_2<\infty$ instead, they reduce to the first point. This is the case, e.g. when a surrogate forward operator $\widetilde{\calG}$ forms a convergent series in the number of training points $n$ (see \citealp{stuart2018posterior}). The IMQ weight in \eqref{eq:w_squared} satisfies the assumptions concerning the weight $w$ with 
\begin{equation*}
    C_s(w)=\frac{1}{\lambda_{\min}}\left(\big\|\widetilde{\calG}(\btheta_0)\big\|^2_2+q^2\lambda_{\max}\right)\ \text{and}\ \gamma_s(\btheta)=\left(\big\|\widetilde{\calG}(\btheta_0)\big\|_2+\frac{q^2}{2}\sqrt{\lambda_{\max}}\right)\frac{\big\|\widetilde{\calG}(\btheta)\big\|^2_2}{\lambda_{\min}}
\end{equation*}
for $\lambda_{\max}$ and $\lambda_{\min}$ the largest and smallest eigenvalues of $\Gamma(\btheta)$ uniformly over $\btheta\in T$.

The stated results ensure perturbation of the posterior by a diverging observational error remains bounded. The cost is paid in the well-specified case, so negligible approximation error and true Gaussian observational noise.

\subsection{Posterior consistency and accuracy}\label{sec:the_acc}

We investigate posterior consistency in the sense of vanishing Hellinger distance of the the regular posterior and posterior approximations for vanishing epistemic approximation error following \citet{stuart2018posterior} and \citet{lie2018random}. The result establishes a trade-off in source of error via the weight $w$. To recover consistency for $\widetilde{\calG}\to\calG$, e.g. in limit of the number of training points $n$, we require $\beta w\to1$. In the non-asymptotic setting, we investigate accuracy between the limits of error divergence and consistency. The upper bounds in Hellinger distance additionally provide upper bounds on the total variation (TV) distance and expected error (see \citealp{dashti2013bayesian}).

We consider an arbitrary, fixed observation $\by\in\mathbb{R}^{d_\by}$ when not stated otherwise. Furthermore, we assume a Gaussian observation likelihood following \eqref{eq:obs_likelihood}. Let $\pi_0$ be the prior measure. We assume $\pi_0$ is absolutely continuous w.r.t to the Lebesgue measure with density $p_0$ considered before. We define the posterior measures under doubling notation via
\begin{itemize}
    \item $\frac{\mathrm{d}p}{\mathrm{d}\pi_0}(\btheta)=\frac{1}{Z}\exp\left(-\frac{1}{2}\|\by-\calG(\btheta)\|^2_{\Gamma_{\boldsymbol\eta}}\right)$ with $Z=\mathbb{E}_{\pi_0}\left[\exp\left(-\frac{1}{2}\|\by-\calG(\btheta)\|^2_{\Gamma_{\boldsymbol\eta}}\right)\right]$;\\
    the regular posterior with true parameter-to-observation map,
    \item $\frac{\mathrm{d}\widetilde{p}}{\mathrm{d}\pi_0}(\btheta)=\frac{1}{\widetilde{Z}}\exp\left(-\frac{1}{2}\|\by-\widetilde{\calG}(\btheta)\|^2_{\Gamma_{\boldsymbol\eta}}\right)$ with $\widetilde{Z}=\mathbb{E}_{\pi_0}\left[\exp\left(-\frac{1}{2}\|\by-\widetilde{\calG}(\btheta)\|^2_{\Gamma_{\boldsymbol\eta}}\right)\right]$;\\ the regular posterior with approximate forward operator, and
    \item $\frac{\mathrm{d}\widetilde{p}_{\beta}}{\mathrm{d}\pi_0}(\btheta)=\frac{1}{\widetilde{Z}_{\beta}}\exp\left(-\frac{1}{2}\beta w^2(\by,\btheta)\|\by-\widetilde{\calG}(\btheta)\|^2_{\Gamma_{\boldsymbol\eta}}\right)$\\
    with $\widetilde{Z}_{\beta}=\mathbb{E}_{\pi_0}\left[\exp\left(-\frac{1}{2}\beta w^2(\by,\btheta)\|\by-\widetilde{\calG}(\btheta)\|^2_{\Gamma_{\boldsymbol\eta}}\right)\right]$;\\ the dynamic WCE posterior with approx. forward operator and learning rate $\beta$.
\end{itemize}
Let $d_\mathrm{H}$ denote the Hellinger distance with
\begin{equation*}
    d^2_\mathrm{H}\left[p(\,\cdot \, | \, \by)\,,\,\widetilde{p}(\,\cdot\,|\by)\right]=\frac{1}{2}\int_T\left(\sqrt{\frac{\mathrm{d}p}{\mathrm{d}\pi_0}}-\sqrt{\frac{\mathrm{d}\widetilde{p}}{\mathrm{d}\pi_0}}\right)^2\pi_0(d\btheta)\,.
\end{equation*}
We additionally denote by 
\vspace{-0.4cm}
\begin{itemize}
    \item $\mathcal{E}_{\boldsymbol{\zeta}}^2(\btheta)\coloneqq\|\calG(\btheta) - \widetilde{\calG}(\btheta)\|_{\Gamma_{\boldsymbol\eta}}^2$ the epistemic $L^2$-error,
    \vspace{-0.2cm}
    \item $\mathcal{E}^2_{\boldsymbol{\eta}}(\btheta,\by)\coloneqq\|\by-\calG(\btheta)\|^2_{\Gamma_{\boldsymbol\eta}}$ the aleatoric $L^2$-error, and 
    \vspace{-0.2cm}
    \item $\mathcal{E}^2(\btheta,\by)\coloneqq\|\by-\widetilde{\calG}(\btheta)\|^2_{\Gamma_{\boldsymbol\eta}}$ the (combined) observation $L^2$-error.
    \vspace{-0.3cm}
\end{itemize}

The normalizing constant $\widetilde{Z}$ is positive and finite under Assumption~\ref{ass:norm_surr} (see \citealp[Lemma 4.1]{stuart2018posterior} for details).
\begin{assumption}\label{ass:norm_surr}
    $\calG$ and $\widetilde{\calG}$ are such that 
    \vspace{-0.4cm}
    \begin{itemize}
        \item $\sup_{\btheta\in T}\|\calG(\btheta)\|^2<\infty$ and
        \vspace{-0.2cm}
        \item $\sup_{\btheta\in T}\mathcal{E}_{\boldsymbol{\zeta}}^2(\btheta)\to0$ in a suitable limit.
        \vspace{-0.3cm}
    \end{itemize}
\end{assumption}
The second assumption may equivalently be stated as $\sup_{\btheta\in T}\mathcal{E}_{\boldsymbol{\zeta}}^2(\btheta)<\infty$ or $\sup_{\btheta\in T}\|\widetilde{\calG}(\btheta)\|^2<\infty$.

Given Assumption~\ref{ass:norm_surr}, \citet[Theorem 4.2]{stuart2018posterior} states
\begin{equation}\label{eq:surr_hell}
    d^2_\mathrm{H}\left[p(\,\cdot \, | \, \by)\,,\,\widetilde{p}(\,\cdot \, | \, \by)\right]\leq C_1\ \left\| \mathcal{E}_{\boldsymbol{\eta}}(\btheta,\by)-\mathcal{E}(\btheta,\by) \right\|^2_{L^2_{\pi_0}(T)} \leq C_2\ \left\| \mathcal{E}_{\boldsymbol{\zeta}}(\btheta)\right\|^2_{L^2_{\pi_0}(T)}
\end{equation}
with $C_1,C_2<\infty$ constants independent of $\btheta\in T$, yet depending on $\by\in\bbR^{d_\by}$. Posterior consistency in $d^2_\mathrm{H}\left[p(\,\cdot \, | \, \by)\,,\,\widetilde{p}(\,\cdot \, | \, \by)\right]\to0$ follows for $\widetilde{\calG}\to\calG$ (or, equivalently, $\mathcal{E}_{\boldsymbol{\zeta}}(\btheta)\to0$), e.g. in the limit of training samples for GP surrogates with \matern covariance kernel and under regularity and boundedness assumptions on $\calG$.

For the normalizing constant of the dynamic WCE posterior, Assumption~\ref{ass: kernel} on $w$ can replace  Assumption~\ref{ass:norm_surr}.

\begin{lemma}[Finite normalizing constant]\label{lem:norm_WCE}
    Let the weight $w$ be such that it satisfies Assumption~\ref{ass: kernel} and $\beta\in\{1,2\}$. The normalizing constant of the dynamic WCE posterior is controlled via
    \begin{equation*}
        \exp(-C(w))\leq \widetilde{Z}_{\beta}\leq1.
    \end{equation*}  
\end{lemma}

\begin{proof}
    Via Assumption~\ref{ass: kernel}, we observe 
    \begin{equation*}
        0\leq \beta w^2(\by,\btheta)\|\by-\widetilde{\calG}(\btheta)\|^2_{\Gamma_{\boldsymbol\eta}}\leq 2C(w)\,.
    \end{equation*}
    Via monotonicity of the exponential function and $\int_T \pi_0(\mathrm{d}\btheta)=1$, the statement follows.
\end{proof}

This enables a corresponding upper bound on the Hellinger distance.

\begin{theorem}[Hellinger distance upper bound]\label{thm: WCE posterior hellinger}
     Let the weight $w$ be such that it satisfies Assumption~\ref{ass: kernel}, $\beta\in\{1,2\}$, and let $\by\in\mathbb{R}^{d_\by}$ be some arbitrary, fixed observation. Then
     \begin{align*}
        d_\mathrm{H}\left[p(\,\cdot \, | \, \by)\,,\,\widetilde{p}_{\beta}(\,\cdot \, | \, \by)\right]&\leq \widetilde{C}_{\mathrm{H}}\ \big\|\mathcal{E}_{\boldsymbol{\eta}}(\btheta,\by)-\sqrt{\beta}w(\by,\btheta)\mathcal{E}(\btheta,\by)\big\|_{L^2_{\pi_0}(T)}\\
        &\leq\widetilde{C}_{\mathrm{H}}\ \left\|\left[1-\sqrt{\beta}w(\by,\btheta)\right]\mathcal{E}_{\boldsymbol{\eta}}(\btheta,\by)\right\|_{L^2_{\pi_0}(T)}\\
          &\hspace{1cm}+\widetilde{C}_{\mathrm{H}}\ \left\|\sqrt{\beta}w(\by,\btheta)\mathcal{E}_{\boldsymbol{\zeta}}(\btheta)\right\|_{L^2_{\pi_0}(T)}
     \end{align*}
     with constant $\widetilde{C}_{\mathrm{H}}=\sqrt{2}\exp\left(\frac{1}{2}(C(w)-1)\right)$.
\end{theorem}

\begin{proof}
    We adjust the proof in \citet{stuart2018posterior}. We observe,
    \begin{align}\label{eq:WCE_hell_split}\nonumber
        d^2_\mathrm{H}\left[p(\,\cdot \, | \, \by)\,,\,\widetilde{p}_{\beta}(\,\cdot \, | \, \by)\right]&\leq\frac{1}{\widetilde{Z}_{\beta}}\int_T \left(\exp\left[-\frac{1}{4}\mathcal{E}^2_{\boldsymbol{\eta}}(\btheta,\by)\right]-\exp\left[-\frac{\beta w^2(\by,\btheta)}{4}\mathcal{E}^2(\btheta,\by)\right]\right)^2\pi_0(d\btheta)\\ 
        &\quad+Z\left(Z^{-\frac{1}{2}}-\widetilde{Z}_{\beta }^{-\frac{1}{2}}\right)^2\eqqcolon(I)+(II)\,.
    \end{align}

    Term $(I)$ is controlled via Lemma~\ref{lem:norm_WCE} and global Lipschitz continuity of $x\mapsto\exp\left(-\frac{1}{4}x^2\right)$ with Lipschitz constant $L=\frac{1}{\sqrt{2e}}$ resulting in
    \begin{equation*}
        (I)\leq\exp(C(w)-1)\big\|\mathcal{E}_{\boldsymbol{\eta}}(\btheta,\by)-\sqrt{\beta}w(\by,\btheta)\mathcal{E}(\btheta,\by)\big\|^2_{L^2_{\pi_0}(T)}
    \end{equation*}
    Regarding term $(II)$, exploit $Z^{\frac{1}{2}}=\left\|\exp\left[-\frac{1}{4}\mathcal{E}^2_{\boldsymbol{\eta}}(\btheta,\by)\right]\right\|_{L^2_{\pi_0}(T)}$ and $\widetilde{Z}_\beta^{\frac{1}{2}}$ accordingly, and via the reverse triangle inequality resulting in
    \begin{equation*}
        (II)\leq\frac{1}{\widetilde{Z}_{\beta}}\left(\widetilde{Z}_{\beta}^{\frac{1}{2}}-Z^{\frac{1}{2}}\right)^2\leq(I)\,.
    \end{equation*}
    
    Placing $(I)$ and $(II)$ in \eqref{eq:WCE_hell_split} gives the stated result, so
    \begin{align*}
        &d^2_\mathrm{H}\left[p(\,\cdot \, | \, \by)\,,\,\widetilde{p}_{\beta}(\,\cdot \, | \, \by)\right]\leq2\cdot(I)\leq2\exp(C(w)-1)\big\|\mathcal{E}_{\boldsymbol{\eta}}(\btheta,\by)-\sqrt{\beta}w(\by,\btheta)\mathcal{E}(\btheta,\by)\big\|^2_{L^2_{\pi_0}(T)}\\
        &\iff d_\mathrm{H}\left[p(\,\cdot \, | \, \by)\,,\,\widetilde{p}_{\beta}(\,\cdot \, | \, \by)\right]\leq\widetilde{C}_{\mathrm{H}}\big\|\mathcal{E}_{\boldsymbol{\eta}}(\btheta,\by)-\sqrt{\beta}w(\by,\btheta)\mathcal{E}(\btheta,\by)\big\|_{L^2_{\pi_0}(T)}
    \end{align*}
    with 
    \begin{equation*}
        \widetilde{C}_{\mathrm{H}}\coloneqq\sqrt{2}\exp\left(\frac{1}{2}(C(w)-1)\right)\,.
    \end{equation*}

     Finally, we observe
    \begin{align*}
        \frac{1}{\widetilde{C}_{\mathrm{H}}}d_\mathrm{H}\left[p(\,\cdot \, | \, \by)\,,\,\widetilde{p}_{\beta}(\,\cdot \, | \, \by)\right]&\leq \left\|\left[1-\sqrt{\beta}w(\by,\btheta)\right]\mathcal{E}_{\boldsymbol{\eta}}(\btheta,\by)\right\|_{L^2_{\pi_0}(T)}\\
        &\quad+\left\|\sqrt{\beta}w(\by,\btheta)\mathcal{E}_{\boldsymbol{\zeta}}(\btheta)\right\|_{L^2_{\pi_0}(T)}\,.
    \end{align*}
    via introducing the term $\sqrt{\beta}w(\by,\btheta)\mathcal{E}_{\boldsymbol{\eta}}(\btheta,\by)-\sqrt{\beta}w(\by,\btheta)\mathcal{E}_{\boldsymbol{\eta}}(\btheta,\by)$, a triangle inequality, and the reverse triangle inequality.
\end{proof}

\begin{remark}\label{rem:const}
    Direct comparison of the statement in Theorem~\ref{thm: WCE posterior hellinger} to \eqref{eq:surr_hell} in \citet[Theorem 4.2]{stuart2018posterior} is limited since Assumption~\ref{ass: kernel} in place of Assumption~\ref{ass:norm_surr} causes the constants to differ in dependencies and scaling. Yet, comparison of the error terms illustrates effect of the weight $w$.
\end{remark}

\begin{corollary}[Finite Hellinger bound]\label{cor:Hell_bound}
    Given the assumptions of Theorem~\ref{thm: WCE posterior hellinger} and let the true forward operator $\calG$ be such that $\sup_{\btheta\in T}\|\calG\|_2<\infty$. Then
    \begin{equation*}
        d_\mathrm{H}\left[p(\,\cdot \, | \, \by)\,,\,\widetilde{p}_{\beta}(\,\cdot \, | \, \by)\right]<\infty\,.
    \end{equation*}
\end{corollary}

\begin{proof}
    We observe
    \begin{align*}
        d^2_\mathrm{H}\left[p(\,\cdot \, | \, \by)\,,\,\widetilde{p}_{\beta}(\,\cdot \, | \, \by)\right]&\leq \widetilde{C}^2_{\mathrm{H}}\ \sup_{\btheta\in T}\left(\mathcal{E}_{\boldsymbol{\eta}}(\btheta,\by)-\sqrt{\beta}w(\by,\btheta)\mathcal{E}(\btheta,\by)\right)^2\\
        &\leq 2\widetilde{C}^2_{\mathrm{H}}\ \max\left\{C_\calG,2C(w)^2\right\}<\infty
    \end{align*}
    via Theorem~\ref{thm: WCE posterior hellinger} with $C_\calG\coloneqq\sup_{\btheta\in T}\mathcal{E}^2_{\boldsymbol{\eta}}(\btheta,\by)\leq\sup_{\btheta\in T}\|\calG(\btheta)\|^2_2+\|\by\|^2_2$ and\\
    $\beta w(\by,\btheta)\mathcal{E}^2(\btheta,\by)\leq 2C(w)$ via Assumption~\ref{ass: kernel}.
\end{proof}

Corollary~\ref{cor:Hell_bound} provides a complement to Theorem~\ref{thm:WCEgbr} in stating robustness in epistemic approximation error. Yet, where \eqref{eq:surr_hell} is consistent in the limit $\mathcal{E}_{\boldsymbol{\zeta}}\to0$, the aleatoric error term prohibits this in Theorem~\ref{thm: WCE posterior hellinger}. Accordingly, consistency of the dynamic WCE posterior additionally requires $\sqrt{\beta}w\to1$.

The statement in Theorem~\ref{thm: WCE posterior hellinger} does not consider stochasticity of the observation. In \eqref{eq:surr_hell} the observation $\by\in\bbR^{d_\by}$ enters the bound via the constant. Assumption~\ref{ass: kernel} allows an expression of the constant independent of $\by$, yet requires explicit control of the stochastic error $\mathcal{E}_{\boldsymbol{\eta}}$ and the weight $w$.

\begin{assumption}\label{ass:err_mom}
    Let the DGD $\pi^\dagger$ be such that
    \begin{equation*}
        m_{\boldsymbol{\eta}}\coloneqq\left\|\bbE_{\pi^\dagger}\left[\mathcal{E}^4_{\boldsymbol{\eta}}(\btheta,\boldsymbol{Y})\right]^{\frac{1}{4}}\right\|_{L^2_{\pi_0}(T)}<\infty\,.
    \end{equation*}
\end{assumption}
For $\sup_{\btheta\in T}\|\calG\|^2_2<\infty$, Assumption~\ref{ass:err_mom} reduces to finite fourth moment of the true observational noise. 

\begin{corollary}[Expected upper bound]\label{cor:exp_bound}
    Given the assumptions of Theorem~\ref{thm: WCE posterior hellinger} and Assumption~\ref{ass:err_mom}. Let $\pi^\dagger$ denote the true DGD. Then
    \begin{align*}
        \bbE_{\pi^\dagger}\left[d_{\rm H}\left[p(\,\cdot\,|\,\boldsymbol{Y}),\widetilde{p}_\beta(\,\cdot\,|\,\boldsymbol{Y})\right]\right]&\leq m_{\boldsymbol{\eta}}\widetilde{C}_{\mathrm{H}}\ \sup_{\btheta\in T}\bbE_{\pi^\dagger}\big[\big(1-\sqrt{\beta}w(\boldsymbol{Y},\btheta)\big)^4\big]^{\frac{1}{4}}\\
        &\quad+\widetilde{C}_{\mathrm{H}}\ \sup_{\btheta\in T}\bbE_{\pi^\dagger}\left[\beta w^2(\boldsymbol{Y},\btheta)\right]^{\frac{1}{2}}\,\|\mathcal{E}_{\boldsymbol{\zeta}}(\btheta)\|_{L^2_{\pi_0}(T)}
        \end{align*}
    for $\widetilde{C}_{\mathrm{H}}<\infty$ the constant in Theorem~\ref{thm: WCE posterior hellinger}.
\end{corollary}

\begin{proof}
    Via Theorem~\ref{thm: WCE posterior hellinger}, linearity of the expectation, Hölder's inequality and Tonelli, we observe
    \begin{align*}
        \bbE_{\pi^\dagger}\left[d_{\rm H}\left[p(\,\cdot\,|\,\boldsymbol{Y}),\widetilde{p}_\beta(\,\cdot\,|\,\boldsymbol{Y})\right]\right]&\leq\widetilde{C}_{\mathrm{H}}\ \big\|\bbE_{\pi^\dagger}\big[\big(1-\sqrt{\beta}w(\boldsymbol{Y},\btheta)\big)^2\mathcal{E}^2_{\boldsymbol{\eta}}(\btheta,\boldsymbol{Y})\big]^{\frac{1}{2}}\big\|_{L^2_{\pi_0}(T)}\\
        &\quad+\widetilde{C}_{\mathrm{H}}\ \left\|\bbE_{\pi^\dagger}\left[\beta w^2(\boldsymbol{Y},\btheta)\right]^{\frac{1}{2}}\mathcal{E}_{\boldsymbol{\zeta}}(\btheta)\right\|_{L^2_{\pi_0}(T)}\eqqcolon(I)+(II)\,.
    \end{align*}

    Term $(I)$ decouples via Cauchy-Schwarz inequality with respect to $\pi^\dagger$ and reduces via Assumption~\ref{ass:err_mom} resulting in
    \begin{align*}
        (I)&\leq\widetilde{C}_{\mathrm{H}}\ \sup_{\btheta\in T}\bbE_{\pi^\dagger}\left[\left(1-\sqrt{\beta}w(\boldsymbol{Y},\btheta)\right)^4\right]^{\frac{1}{4}}\ \left\|\bbE_{\pi^\dagger}\left[\mathcal{E}^4_{\boldsymbol{\eta}}(\btheta,\boldsymbol{Y})\right]^{\frac{1}{4}}\right\|_{L^2_{\pi_0}(T)}\\
        &\leq m_{\boldsymbol{\eta}}\widetilde{C}_{\mathrm{H}}\ \sup_{\btheta\in T}\bbE_{\pi^\dagger}\left[\left(1-\sqrt{\beta}w(\boldsymbol{Y},\btheta)\right)^4\right]^{\frac{1}{4}}\,.
    \end{align*}

    Term $(II)$ reduces similarly via
    \begin{equation*}
        (II)\leq\widetilde{C}_{\mathrm{H}}\ \sup_{\btheta\in T}\bbE_{\pi^\dagger}\left[\beta w^2(\boldsymbol{Y},\btheta)\right]^{\frac{1}{2}}\left\|\mathcal{E}_{\boldsymbol{\zeta}}(\btheta)\right\|_{L^2_{\pi_0}(T)}\,. 
    \end{equation*}
    The stated result directly follows.
\end{proof}

Consistency under vanishing approximation error $\mathcal{E}_{\boldsymbol{\zeta}}(\btheta)\to0$ requires $\beta w^2\to 1$ to not pick up on the stochastic error instead. Accuracy of the dynamic WCE posterior depends on the ability of the weight $w$ to determine the smaller instance of error to maintain. Theorem~\ref{thm: WCE posterior hellinger} and Corollary~\ref{cor:exp_bound} isolate this trade-off, i.e., for instances of large approximation error the stochastic error of reduced observation assimilation should be maintained. Additionally, the results illustrate the role of the learning rate $\beta$ as it determines whether the condition  $\beta w^2\to1$ reflects the limit case $w\to1$ for $\beta=1$ or an intermediate case $w^2\approx\frac{1}{2}$ for $\beta=2$ which allows to recover consistency while maintaining robustness.

\begin{remark}
    The limit case $w\to1$ for all $(\by,\btheta)\in\bbR^{d_\by}\times T$ requires $q^2\to\infty$ for the IMQ weight in \eqref{eq:w_squared}. This causes the lower bound of the normalizing constant $\widetilde{Z}_\beta$ to vanish, yet it can be recovered via Assumption~\ref{ass:norm_surr}. 
\end{remark}

\paragraph{Transfer to design choices.}

For considering the static parametrisation, results transfer when additionally assuming the conditions on the weight in Corollary~\ref{cor:statgbr} next to Assumption~\ref{ass: kernel} and Assumption~\ref{ass:norm_surr} to maintain boundedness arguments.

For considering surrogates with parameter dependent covariance matrices such as the marginal GP surrogate in \eqref{eq:marginal_surrogate}, the introduced error terms consider different norms in $\|\cdot\|^2_{\Gamma(\btheta)}$ and $\|\cdot\|^2_{\Gamma_{\boldsymbol{\eta}}}$. Lemma~\ref{lem:norm_WCE} translates, yet the statement in Theorem~\ref{thm: WCE posterior hellinger} does not. We state an adjusted result for the dynamic WCE posterior with marginal GP surrogate in \eqref{eq:marginal_surrogate} utilizing $\Gamma(\btheta)=\Gamma_{\boldsymbol{\eta}}+K_n^{\calG}(\btheta,\btheta)$.

\begin{corollary}[Parameter dependent covariance]\label{cor:marg_hell}
    Given the assumptions of Theorem~\ref{thm: WCE posterior hellinger} with the marginal GP surrogate in \eqref{eq:marginal_surrogate} and parameter dependent covariance matrix $\Gamma(\btheta)=\Gamma_{\boldsymbol{\eta}}+K_n^{\calG}(\btheta,\btheta)$. Then
    \begin{align*}
        d_\mathrm{H}\left[p(\,\cdot \, | \, \by),\widetilde{p}_{\beta}(\,\cdot \, | \, \by)\right]&\leq \widetilde{C}_{\mathrm{H}}\ \left\|\left[1-\sqrt{\beta}w(\by,\btheta)\right]\mathcal{E}_{\boldsymbol{\eta}}(\btheta,\by)\right\|_{L^2_{\pi_0}(T)}\\
        &\quad+\widetilde{C}_{\mathrm{H}}\ \left\|\sqrt{\beta}w(\by,\btheta)\mathcal{E}_{\boldsymbol{\zeta}}(\btheta)\right\|_{L^2_{\pi_0}(T)}+\widetilde{C}_{\mathrm{H}}\ \left\|\sqrt{\beta C(w)\widetilde{\lambda}(\btheta)}\right\|_{L^2_{\pi_0}(T)}
     \end{align*}
     for $\widetilde{C}_{\mathrm{H}}<\infty$ the constant in Theorem~\ref{thm: WCE posterior hellinger} and $\widetilde{\lambda}(\btheta)=\lambda_{\max}\left(\Gamma_{\boldsymbol{\eta}}^{\frac{1}{2}}K_n^{\calG}(\btheta,\btheta)\Gamma_{\boldsymbol{\eta}}^{\frac{1}{2}}\right)$. 
\end{corollary}

\begin{proof}
    Let $\widetilde{\mathcal{E}}(\btheta,\by)\coloneqq\|\by-\widetilde{\calG}(\btheta)\|_{\Gamma(\btheta)}$. Theorem~\ref{thm: WCE posterior hellinger} holds in the first inequality with the adjusted error term, so
    \begin{equation*}
        d_\mathrm{H}\left[p(\,\cdot \, | \, \by),\widetilde{p}_{\beta}(\,\cdot \, | \, \by)\right]\leq \widetilde{C}_{\mathrm{H}}\ \left\|\mathcal{E}_{\boldsymbol{\eta}}(\btheta,\by)-\sqrt{\beta}w(\by,\btheta)\widetilde{\mathcal{E}}(\btheta,\by)\right\|_{L^2_{\pi_0}(T)}\,.
    \end{equation*}
    
    We approach the second inequality the same way via
    \begin{align*}
        \frac{1}{\widetilde{C}_{\mathrm{H}}}d_\mathrm{H}\left[p(\,\cdot \, | \, \by),\widetilde{p}_{\beta}(\,\cdot \, | \, \by)\right]&\leq \left\|\left[1-\sqrt{\beta}w(\by,\btheta)\right]\mathcal{E}_{\boldsymbol{\eta}}(\btheta,\by)\right\|_{L^2_{\pi_0}(T)}\\
        &\quad+\left\|\sqrt{\beta}w(\by,\btheta)\left[\mathcal{E}_{\boldsymbol{\eta}}(\btheta,\by)-\widetilde{\mathcal{E}}(\btheta,\by)\right]\right\|_{L^2_{\pi_0}(T)}\,.
    \end{align*}

    We observe
    \begin{equation*}
        \left|\mathcal{E}_{\boldsymbol{\eta}}(\btheta,\by)-\widetilde{\mathcal{E}}(\btheta,\by)\right|\leq\mathcal{E}_{\boldsymbol{\zeta}}(\btheta)+\left\langle\by-\widetilde{\calG}(\btheta),\left(\Gamma_{\boldsymbol{\eta}}^{-1}-\Gamma(\btheta)^{-1}\right)(\by-\widetilde{\calG}(\btheta))\right\rangle^{\frac{1}{2}}
    \end{equation*}
    via introducing the term $\|\by-\widetilde{\calG}(\btheta)\|_{\Gamma_{\boldsymbol\eta}}-\|\by-\widetilde{\calG}(\btheta)\|_{\Gamma_{\boldsymbol\eta}}$, the reverse triangle inequality and the sub-additivity of the square root.

    Let $\by-\widetilde{\calG}(\btheta)=\Gamma(\btheta)^{\frac{1}{2}}\boldsymbol{z}$ and $\widetilde{\mathcal{E}}(\btheta,\by)=\boldsymbol{z}^T\boldsymbol{z}$. Then
    \begin{align*}
        \left\langle\by-\widetilde{\calG}(\btheta),\left(\Gamma_{\boldsymbol{\eta}}^{-1}-\Gamma(\btheta)^{-1}\right)(\by-\widetilde{\calG}(\btheta))\right\rangle&=\left\langle\boldsymbol{z},\left(\Gamma(\btheta)^{\frac{1}{2}}\Gamma_{\boldsymbol{\eta}}^{-1}\Gamma(\btheta)^{\frac{1}{2}}-\mathrm{Id}_{d_\by\times d_\by}\right)\boldsymbol{z})\right\rangle\\
        &\leq\widetilde{\lambda}(\btheta)\widetilde{\mathcal{E}}(\btheta,\by)
    \end{align*}
    with $\widetilde{\lambda}(\btheta)$ the largest eigenvalue of $\Gamma_{\boldsymbol{\eta}}^{\frac{1}{2}}K_n^{\calG}(\btheta,\btheta)\Gamma_{\boldsymbol{\eta}}^{\frac{1}{2}}$ via
    \begin{equation*}
        \lambda_{\max}\left(\Gamma(\btheta)^{\frac{1}{2}}\Gamma_{\boldsymbol{\eta}}^{-1}\Gamma(\btheta)^{\frac{1}{2}}-\mathrm{Id}_{d_\by\times d_\by}\right)=\lambda_{\max}\left(\Gamma_{\boldsymbol{\eta}}^{\frac{1}{2}}K_n^{\calG}(\btheta,\btheta)\Gamma_{\boldsymbol{\eta}}^{\frac{1}{2}}\right)\,.
    \end{equation*}

    Taking both steps together results in 
    \begin{align*}
        d_\mathrm{H}\left[p(\,\cdot \, | \, \by),\widetilde{p}_{\beta}(\,\cdot \, | \, \by)\right]&\leq \widetilde{C}_{\mathrm{H}}\ \left\|\left[1-\sqrt{\beta}w(\by,\btheta)\right]\mathcal{E}_{\boldsymbol{\eta}}(\btheta,\by)\right\|_{L^2_{\pi_0}(T)}\\
        &\quad+\widetilde{C}_{\mathrm{H}}\ \left\|\sqrt{\beta}w(\by,\btheta)\mathcal{E}_{\boldsymbol{\zeta}}(\btheta)\right\|_{L^2_{\pi_0}(T)}+\widetilde{C}_{\mathrm{H}}\ \left\|\sqrt{\beta C(w)\widetilde{\lambda}(\btheta)}\right\|_{L^2_{\pi_0}(T)}
    \end{align*}
    via Assumption~\ref{ass: kernel}.
\end{proof}

\begin{remark}
    The additional term $\widetilde{C}_{\mathrm{H}}\ \big\|\sqrt{\beta C(w)\widetilde{\lambda}(\btheta)}\big\|_{L^2_{\pi_0}(T)}$ is a rough upper bound on the effect of the additive covariance inflation of the marginal GP surrogate independent of the observation $\by\in\bbR^{d_\by}$ via Assumption~\ref{ass: kernel}. Accordingly, taking the expectation w.r.t. the true DGD $\pi^\dagger$ does not affect it and Corollary~\ref{cor:exp_bound} applies up to the additive term. 
\end{remark}

The previous insights in consistency hold in vanishing approximation error via $\widetilde{\lambda}(\btheta)\to0$, e.g. in the limit $n\to\infty$ for the marginal GP surrogate in \eqref{eq:marginal_surrogate} (see \citealp{stuart2018posterior}). Again, the bound is not suited to consider the case $q^2\to\infty$ for the IMQ weight in \eqref{eq:w_squared} to recover uniform $w^2\to1$ since $C(w)\to\infty$, yet highlights the role of the learning rate $\beta$.

\section{Posterior evaluation}
\label{sec:sampling_methodology}
For evaluating the WCE posteriors with GP surrogates we can exploit closed form access to posterior score functions. We consider the Metropolis-adjusted Langevin algorithm (MALA)  \cite{robert1999monte,roberts1996exponential} presented in Algorithm~\ref{alg:mala} as an example.

\begin{algorithm}[tb]
    \DontPrintSemicolon
    Given a (generalised) posterior density $p_{\beta\mathrm{D}}(\,\cdot \,|\, \by)$, choose initial state $\btheta_0$ and step size $\alpha$. \;
    We define $m$ as
    \begin{align*}
        m(\btheta \,|\, \btheta')
        \propto \exp\left(- \frac{1}{4\alpha} \, \left\|\btheta - \btheta' - \alpha \nabla_{\btheta} \log p_{\beta\mathrm{D}}(\btheta' \,|\, \by)\right\|_2^2\right)\,.
    \end{align*}
    
    \For{$j=0,~1,~\ldots,~J$}{
        Draw $\bxi_j$ as independent realisation from $\calN(0, \mathrm{Id}_{d_{\btheta}})$. \;
        Propose $\widehat{\btheta}_j = \btheta_j + \alpha \, \nabla_{\btheta} \log p_{\beta\mathrm{D}}(\btheta_j \,|\, \by) + \sqrt{2 \alpha} \, \bxi_j$. \;
        Set $\btheta_{j+1} = \widehat{\btheta}_j$ with probability $\nu_j$ given by
        \begin{align*}
            \nu_j
            &= \min \left\{1, \frac{p_{\beta\mathrm{D}}(\widehat{\btheta}_j \,|\, \by)}{p_{\beta\mathrm{D}}(\btheta_j \,|\, \by)} \, \frac{m(\btheta_j \,|\, \widehat{\btheta}_j)}{m(\widehat{\btheta}_j \,|\, \btheta_j)}\right\}\,.
        \end{align*}
        
        Otherwise, set $\btheta_{j+1} = \btheta_j$. \;
    }
    \caption{Metropolis-adjusted Langevin algorithm.}
    \label{alg:mala}
\end{algorithm}

\begin{remark}\label{rem:init_state}
    We propose choosing the initial state $\btheta_0$ as the approximate mode of the posterior with the mean GP likelihood in \eqref{eq:mean_surrogate}, e.g. via Adam (see \citealp{kingma2015iclr-adam}), to start MALA closer to its target distribution since the gradients of the WCE posteriors vanish away from the mode.
\end{remark}

\paragraph{Weighted cross-entropy posterior gradients.}
We define
\begin{equation*}
    \mathrm{QF}(\by,\btheta) \coloneqq \|\by - \widetilde{\calG}(\btheta)\|^2_{\Gamma(\btheta)}
\end{equation*}
with $\Gamma(\btheta)$ depending on the choice of GP surrogate and $\widetilde{\ell}$ the corresponding observation likelihood specified by either $\widetilde{\ell}_\mathrm{mean}$ in \eqref{eq:mean_surrogate} or $\widetilde{\ell}_\mathrm{GP}$ in \eqref{eq:marginal_surrogate}.

\begin{proposition}
    \label{prop:wce}
    The score functions of the WCE posteriors are given by
    \begin{align*}
        \nabla_{\btheta} \log \widetilde{p}_{\beta}(\btheta \,|\, \by)
        &= \beta w^2(\by, \btheta) \nabla_{\btheta} \log \widetilde{\ell}_+(\by \,|\, \btheta)
        + \nabla_{\btheta} \log p_0(\btheta)\\
        &\quad+ \beta \log \widetilde{\ell}_+(\by \,|\, \btheta) \nabla_{\btheta} w^2(\by, \btheta)\\
        &\text{and}\\
        \nabla_{\btheta} \log \bar p_{\beta}(\btheta \,|\, \by)
        &= \beta w^2(\by, \bar\btheta) \nabla_{\btheta} \log \widetilde{\ell}_+(\by \,|\, \btheta)
        + \nabla_{\btheta} \log p_0(\btheta)\,.
    \end{align*}
    For $\widetilde{\ell}$ in \eqref{eq:mean_surrogate} or \eqref{eq:marginal_surrogate} and the IMQ weight $w$ in \eqref{eq:w_squared}, the score function of the WCE posteriors are given by
    \begin{align*}
        \nabla_{\btheta} \log \widetilde{p}_{\beta}(\btheta \,|\, \by)&=-\frac{\beta}{2} w^4(\by, \btheta)\nabla_{\btheta} \mathrm{QF}(\by, \btheta) + \nabla_{\btheta} \log p_0(\btheta)\\
        &\text{and}\\
        \nabla_{\btheta} \log \bar p_{\beta}(\btheta \,|\, \by)
        &= -\frac{\beta}{2} w^2(\by, \bar\btheta) \nabla_{\btheta} \mathrm{QF}(\by, \btheta)+ \nabla_{\btheta} \log p_0(\btheta)
    \end{align*}
    with $\nabla_{\btheta}\mathrm{QF}(\by,\btheta)$ available in closed form as in \citet[Lemma 1]{bai2024gaussian}.
\end{proposition}

\begin{proof}
    The first statement follows directly from linearity of the gradient with the product rule in the first equation and $w$ a scalar w.r.t. $\btheta\in T$ in the second equation. The second statement follows via expressing the IMQ weight $w$ in \eqref{eq:w_squared} via  
    \begin{equation*}
        w^2(\by, \btheta)=\left(1 + \frac{\mathrm{QF}(\by,\btheta)}{q^2}\right)^{-1}\quad\text{with}\quad\nabla_{\btheta} w^2(\by,\btheta)= - \frac{w^4(\by,\btheta)}{q^2}\nabla_{\btheta}\mathrm{QF}(\by,\btheta)
    \end{equation*}
    and the identity
    \begin{equation*}
        \left(w^2(\by, \btheta)-\frac{1}{q^2}w^4(\by, \btheta)\mathrm{QF}(\by,\btheta)\right)= w^4(\by, \btheta)\,.
    \end{equation*}
\end{proof}

The second statement in Proposition~\ref{prop:wce} in addition to availability of the GP surrogate gradients in \citet{bai2024gaussian} motivates tractability of the WCE posterior for the presented approaches.

\paragraph{Parameter tuning.}
Most choices of weight $w$ in WCE via kernels introduce an additional parameter as a degree of freedom in exchange for weakened assumptions on the observation likelihood, i.e. not considering it well-specified. For the IMQ weight the parameter $0<q^2<\infty$ guarantees robustness, yet strongly effects consistency and accuracy. Tuning the parameter $q^2$ is challenging and usually discussed in specific contexts only (see \citealp{altamirano2023robustBOCPD,duran2024outlier,reimann2026kalman})

Following the  heuristic in \citet{altamirano2023robustBOCPD}, we may choose $q^2$ according to the WCE posterior with IMQ weight to perform empirically accurate in the well-specified case in the following sense: Given a computational budget for offline tuning, we can generate observations $\{\by_i\}_{i=1}^{M_{\mathrm{train}}}$ via the true parameter-to-observation map $\calG$, with $M_{\mathrm{train}}\in\mathbb{N}$ the number of artificial observations. For a given $0<q^2<\infty$, let $\{\{\btheta_i^{j,q}\}_{j=1}^{M_{\mathrm{MC}}}\}_{i=1}^{M_{\mathrm{train}}}$ be the corresponding samples of $\widetilde{p}_{\beta;q^2}(\,\cdot \,|\, \by_i)$ with sample size $M_{\mathrm{MC}}\in\mathbb{N}$. We define
\begin{equation*}
    \mathcal{L}_\mathrm{L}(q)=\frac{1}{M_\mathrm{train}}\sum_{i=1}^{M_\mathrm{train}}\widehat{\mathrm{L}}\left[\{\btheta_i^{j,q}\}_{j=1}^{M_\mathrm{MC}} \,\big\|\, \by_i\right]
\end{equation*}
with $\mathrm{L}$ some loss function with empirical estimator $\widehat{\mathrm{L}}$ between sample and observation. The posterior cross-entropy $-\log \widetilde{p}_{\beta;q^2}(\btheta \,|\, \by_i)$ (see \citealp{gneiting2007strictly}) is a suitable candidate, yet requires posterior density estimation. Alternatively, $\mathrm{L}$ can be chosen as the estimated KL divergence $\widehat{\mathrm{KL}}\left[p(\btheta \,|\,\by_i) \,\|\, \widetilde{p}_{\beta;q^2}(\btheta \,|\, \by_i)\right]$, e.g. via \cite{perez2008kullback}, which additionally requires access to a sample $\{\{\btheta_i^j\}_{j=1}^{M_{\mathrm{MC}}}\}_{i=1}^{M_{\mathrm{train}}}$ of $p(\, \cdot \,|\,\by_i)$. Parameter choice then follows via $q_\mathrm{L}^*\coloneqq\arg\min_{q^2>0}\mathcal{L}_\mathrm{L}(q)$, e.g. over a predetermined sequence.

From an analytical perspective we may assume \eqref{eq:obs_likelihood} conditioned on $\btheta$ generating $\boldsymbol{Y}$ is well-specified. Then $-2\log\ell_+(\boldsymbol{Y}\,|\,\btheta)\sim\chi^2(d_\by)$ and $\bbE_{\boldsymbol{Y}|\btheta}\left[-\log\ell_+(\boldsymbol{Y}\,|\,\btheta)\right]=\frac{1}{2}d_{\by}$. We observe, $\|\widetilde{\calG}(\btheta)-\boldsymbol{Y}\|^2_{\Gamma_{\boldsymbol{\eta}}}\sim\chi^2(d_\by;\mathcal{E}_{\boldsymbol{\zeta}}^2(\btheta))$, the non-central $\chi^2$-distribution with $d_\by$ degrees of freedom and non-centrality parameter $\mathcal{E}_{\boldsymbol{\zeta}}^2(\btheta)$. Similarly, we observe $\|\widetilde{\calG}(\btheta)-\boldsymbol{Y}\|^2_{\Gamma(\btheta)}\sim\chi^2(d_\by)$ for assuming $\calG(\btheta)-\widetilde{\calG}(\btheta)\sim\mathcal{N}(0,\Gamma(\btheta)-\Gamma_{\boldsymbol{\eta}})$. We adjust \citet[Lemma 1]{reimann2026kalman} as a special case of \citet{liao2018sharpening} and \citet{gao2017bounds}.

\begin{lemma}\label{lem:tuning}
    Suppose the non-negative random variable $Z$ satisfies $\bbE[Z]=\mu\geq0$ and $\mathrm{Var}(Z)=\sigma^2<\infty$. Let $h:[0,\infty)\to\bbR;\ z\mapsto z\left(1+\frac{z}{q^2}\right)^{-1}$ for $0<q^2<\infty$. Then
    \begin{equation*}
        \left(1+\frac{\mu}{q^2}\right)^{-1}\mu-\frac{q^2\sigma^2}{(q^2+\mu)^2}\leq\bbE[h(Z)]\leq\left(1+\frac{\mu}{q^2}\right)^{-1}\mu\,.
    \end{equation*}
\end{lemma}

\begin{proof}
    Since $h$ is concave on $[0,\infty)$ and via Jensen's inequality, $\bbE[h(Z)]\leq h(\mu)$. For the lower bound, we observe
    \begin{align*}
       h(\mu)-\bbE\left[h(Z)\right]&=\bbE\left[h(\mu)-h(Z)-h'(\mu)(\mu-Z)\right]=q^4\bbE\left[\frac{(\mu-Z)}{(q^2+\mu)(q^2+Z)}-\frac{(\mu-Z)}{(q^2+\mu)^2}\right]\\
       &=q^4\bbE\left[\frac{(Z-\mu)^2}{(q^2+\mu)^2(q^2+Z)}\right]\leq\frac{q^2\sigma^2}{(q^2+\mu)^2}
    \end{align*}
    since $\bbE[\mu-Z]=0$, $h'(z)=\left(1+\frac{z}{q^2}\right)^{-2}$, and $q^2+Z\geq q^2$.
\end{proof}
We define $Z\coloneqq-2\log\widetilde{\ell}_+(\boldsymbol{Y}\,|\,\btheta)=\mathrm{QF}(\boldsymbol{Y},\btheta)$. Then $-\beta w^2(\by,\btheta)\log\widetilde{\ell}_+(\by\,|\,\btheta)=h(Z)$ with $h$ as in Lemma~\ref{lem:tuning} for $\beta=2$. Matching the upper bound in Lemma~\ref{lem:tuning} to recover the expected rescaled true negative log-likelihood $\frac{1}{2}d_\by$, we observe for the approximate likelihood without parameter dependent covariance
\begin{equation*}
    \frac{1}{2}d_\by-\frac{d_\by[d_\by+2\mathcal{E}^2_{\boldsymbol{\zeta}}(\btheta)]^2}{2[d_\by+\mathcal{E}_{\boldsymbol{\zeta}}^2(\btheta)]^3}\leq\bbE[h(Z)]\leq\frac{1}{2}d_\by\quad\text{for}\quad q^2=\frac{d_\by\left[d_\by+\mathcal{E}_{\boldsymbol{\zeta}}^2(\btheta)\right]}{d_\by+2\mathcal{E}_{\boldsymbol{\zeta}}^2(\btheta)}
\end{equation*}
and for the approximate likelihood with parameter dependent covariance 
\begin{equation*}
    \frac{1}{2}d_\by-\frac{1}{2}\leq\bbE[h(Z)]\leq\frac{1}{2}d_\by\quad\text{for}\quad q^2=d_\by\,.
\end{equation*}
The approximation error $\mathcal{E}_{\boldsymbol{\zeta}}^2(\btheta)$ is generally unknown as motivates the presented approaches. For the GP mean surrogate, we propose estimation via the GP covariance kernel, so $\widehat{\mathcal{E}}^2_{\boldsymbol{\zeta}}=\frac{1}{M_{\mathcal{E}}}\sum_{i=1}^{M_{\mathcal{E}}}\mathrm{Tr}(\Gamma_{\boldsymbol\eta}^{-1}K_n^{\calG}(\btheta_i,\btheta_i))$ for a sample $\{\btheta_i\}_{i=1}^{M_{\mathcal{E}}}$ of the prior $p_0$ with $M_{\mathcal{E}}\in\mathbb{N}$. 

The analytic approach is an idealized setting and strongly depends on the underlying assumptions. It does not cover the static WCE or a learning rate $\beta=1$. The proposed values provide a justified estimate without rigorous guarantees. However, we can conjecture that $q^2$ needs to scale with the observation dimension $d_\by$ which can inform the grid in offline tuning.  

\section{Numerical experiments}
\label{sec:numerical_experiment}

We provide a basic study of design choices and consistency to complement theory and guide further considerations. Since the weighted cross-entropy posterior aims at improving UQ, we focus evaluation on metrics capturing sample spread around a true parameter.

We follow the setup described in \citet[Section 4.1.1]{bai2024gaussian} and consider a one-dimensional elliptic PDE
\begin{equation}
    \label{eq:elliptic_pde}
    \begin{split}
        - \frac{\mathrm{d}}{\mathrm{d}x} \left(\exp(\kappa(x, \btheta)) \, \frac{\mathrm{d}}{\mathrm{d}x} u(x)\right)
        &= 4 x
        \quad \text{for }\quad x \in (0, 1)\,, \\
        u(0) &= 0\,, \quad \text{and} \quad u(1) = 2
    \end{split}
\end{equation}
with function $\kappa$ governing the diffusion coefficient defined as piece-wise constant via
\begin{equation*}
    \kappa(x, \btheta) = \vartheta_1\,\mathbf{1}_{[\frac{1}{4},\,\frac{1}{2})}(x) + \vartheta_2\,\mathbf{1}_{[\frac{1}{2},\,\frac{3}{4})}(x) + \mathbf{1}_{[\frac{3}{4},\,1)}(x)
    \quad\text{for}\quad\btheta=(\vartheta_1,\vartheta_2)^T\in T=[-1, 1]^2\,.
\end{equation*}

The prior given is by a uniform distribution on $[-1, 1]^2$, which when applying Algorithm \ref{alg:mala} is smoothed by convolution with a Gaussian kernel to ensure $\nabla_{\btheta}\log p_0$ is well-defined. The true parameter $\btheta^\dagger$ is drawn from the prior.

Artificial observations $\by$ are generated from the solution $u$ given $\btheta^\dagger$ via a finite element solver on $2560$ elements observed at $d_{\by}=6$ equally-spaced points in $[0,1]$ and adding Gaussian observational noise with $\Gamma_{\boldsymbol\eta} = 10^{-6} \, \mathrm{Id}_{d_{\by}}$.

\begin{figure}[tb]
    \centering
    \includegraphics[width=\linewidth]{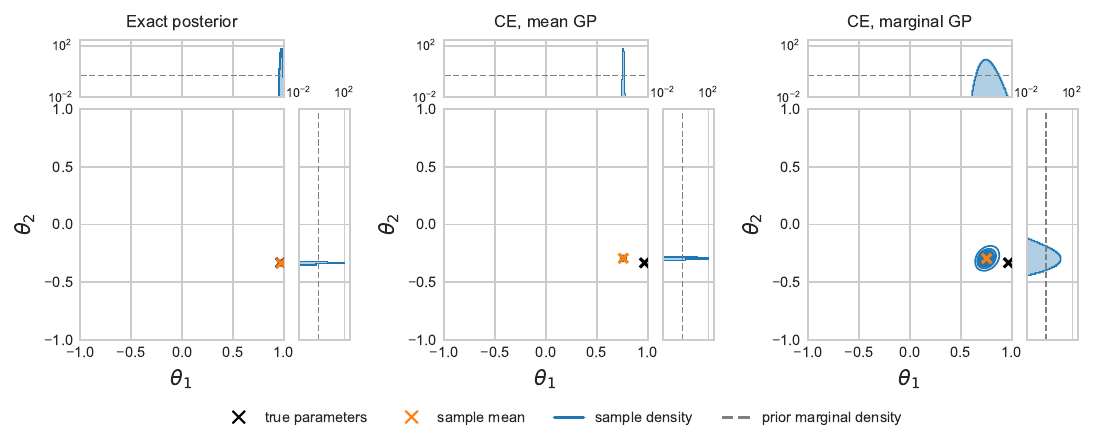}
    \caption{Joint posterior density contours of $\btheta$ with marginal histograms (log-frequency) of regular posterior samples with the true forward operator $\calG$ (left), the mean GP surrogate (middle), and the marginal GP surrogate (right).}
    \label{fig:ce_exact_surr}
\end{figure}

\paragraph{Regular posterior evaluation.}
We compute samples of the posterior in \eqref{eq:posterior} employing the true forward operator $\calG$ based on solutions of \eqref{eq:elliptic_pde} via a finite element solver on $320$ elements and running a Metropolis-Hastings algorithm \citep{hastings1970monte,metropolis1953equation} for $10^6$ steps and discard the first half of the obtained samples as a burn-in phase. A density plot of the posterior samples obtained this way is given in Figure~\ref{fig:ce_exact_surr}. This regular posterior provides the most accurate approximation to the true posterior assumed available, avoids an inverse crime, and serves to indicate difficulty of the inverse problem.

\paragraph{GP surrogate posterior evaluation.}
We replace solving \eqref{eq:elliptic_pde} numerically at each step of our sampling algorithm with a GP surrogate approximation. The GP surrogate is pre-trained for $n\in\mathbb{N}$ training points in $[-1, 1]^2$ using the Halton sequence with $\calG$ based on solutions of \eqref{eq:elliptic_pde} via a finite element solver on $320$ elements. The GP surrogate model uses a squared exponential covariance kernel and the hyper-parameters are selected via maximising the marginal likelihood. We compute samples of the all approximate posteriors via MALA with the gradient expressions given by \cite{bai2024gaussian} or in Proposition~\ref{prop:wce}, for $10^6$ steps and discard the first half of the samples as burn-in. The initial values are chosen following Remark~\ref{rem:init_state} with Adam initialized at the prior mean.  A density plot of the posterior samples obtained this way is given in Figure~\ref{fig:ce_exact_surr}. We compare across configurations of the introduced modular framework as well as with the regular posterior with either GP likelihood approximation.

\paragraph{Evaluation metrics.}
We evaluate uncertainty quantification via a nearest-neighbour estimator of the negative log posterior density (NLPD) at the true parameter $\btheta^\dagger$ based on \citet{loftsgaarden1965nonparametric} with $k=\lceil\sqrt{M}\rceil$-many neighbours and the empirical coverage of the marginal $90\%$ credible intervals. All metrics are averaged over $M_{\mathrm{rep}}=100$ repetitions of each simulation setting. Next to the regular posterior with the FEM solve of the forward operator, we provide all metrics with respect to the prior to provide a spectrum of values for the NLPD.

\paragraph{Comparison of design choice.}
We employ $n=5$ training points for the GP surrogates since that covers a middle ground where the surrogates are useful yet exhibit impactful approximation error. Implementation follows Section~\ref{sec:sampling_methodology} with both tuning approaches. Offline tuning uses $M_{\mathrm{MC}}=10^5$ samples and $M_{\mathrm{train}}=50$ repetitions for each grid point with $\widetilde{q}\in\{0.5, 0.75, 1, 2, 3, 5, 10, 20, 50, 100\}$ and $q^2=\widetilde{q}^2d_\by$. The grid is chosen intentionally coarse to enable comparison across different number of training points $n$. The values $\widetilde{q}=0.5$ and $\widetilde{q}=100$ essentially recover the prior and regular CE posterior in practice. The artificial observations $\by_i$ are generated via $\btheta_i$ drawn from the prior and employed in the forward model used for the GP training data, so $\calG$ based on solutions of \eqref{eq:elliptic_pde} via a finite element solver on $320$ elements. The tuned $q^2$ minimizes the average empirical NLPD and values are reported in Table~\ref{tab:grid_q}. The heuristic tuning via Lemma~\ref{lem:tuning} suggests $q^2=3.01$ for the mean GP surrogate with approximation error $\widehat{\mathcal{E}}^2_{\boldsymbol{\zeta}}=0.708$ using $M_{\mathcal{E}}=10^4$ prior samples and $q^2=d_\by=6$ for the marginal GP surrogate.

\input{results/grid_q}

The results are presented in Table~\ref{tab:results_NLL_90}. The static weighted cross-entropy posteriors with $q^2$ tuned via the grid search reliably improve inference for both surrogates and learning rates. Its strategy via observation driven scaling of the observation likelihood covariance provides the lowest NLPD at the true parameter (which is also the grid-tuning criterion) and improved coverage except for the marginal GP surrogate with $\beta=2$. The dynamic WCE posteriors tends to recover the prior density at the true parameter for the mean GP surrogate with slight improvement over the prior for the marginal GP surrogate with overconfidence in coverage. Different from the static weighting, the dynamic weights strongly alter the shape of the posterior. The results suggest that choosing $w$ as the IMQ weight in \eqref{eq:w_squared} does not suit the dynamic weighting and alternative specifications need to be considered. Finally, the heuristic, or proposed default, choice of $q^2$ does generally not lead to better inference over the regular CE posteriors to a relevant degree apart from inconsistent improvements for the static weighting with conservative coverage. Considering that the assumption on which it is based strongly limits the heuristic tuning in applicability in the first place, this is within expectation. The weighted posterior and the tuning heuristic are points of discussion next.

\input{results/table_NLL_90}


\section{Discussion, challenges and conclusion}\label{sec:disc_conclusion}

The proposed approach cannot change the bias induced by the approximation error directly. With the uniform prior, the weighted cross-entropy posterior maintains the cross-entropy posterior mode up to the log-determinate factor for the marginal GP surrogate (RMSE of $0.094$ for every mean GP surrogate and $0.093$ to $0.096$ for the marginal GP surrogates). Note however, that employing WCE affects the posterior mean and pulls it towards the prior mean which shifts it from the CE posterior and can coincidentally improve mean estimation (RMSE of $0.093$ to $0.094$ for the CE posteriors and $0.078$ to $0.093$ for the static WCE posteriors with grid tuning). In the context of standard Bayesian inference, improving on the accuracy of the posterior mode was achieved by employing better prior distributions for the GP surrogate model in \citet{bai2024gaussian}, and we expect to see similar improvements for generalized Bayesian inference.

By acknowledging the presence of the parameter-dependent bias as an additional source of error, we can improve uncertainty quantification and reduce overconfidence. The weighted cross-entropy posteriors bridge between the prior and the regular cross-entropy posterior with the surrogate forward operator. The former is unaffected by the approximation error but uninformative. The latter is informative but highly overconfident in the biased mode (see Table~\ref{tab:results_NLL_90}; $4\%$ and $38\%$ coverage at the $90\%$ level). This recovers a from of bias-variance trade-off. The presented approach addresses the resulting challenge to balance both terms as we observe in theory in Section~\ref{sec:the_acc} and application in Section~\ref{sec:numerical_experiment}. Utilizing the weighted cross-entropy in \citet{duran2024outlier} in the context of generalized Bayesian inference shifts this challenge to specification of the weight $w$, and for the IMQ weight in \eqref{eq:w_squared} it shifts again to tuning the parameter $q^2$. We emphasize that the presented components are deliberate design choices and can likely be improved in different contexts. 

\begin{figure}[tb]
    \centering
    \includegraphics[width=0.67\linewidth]{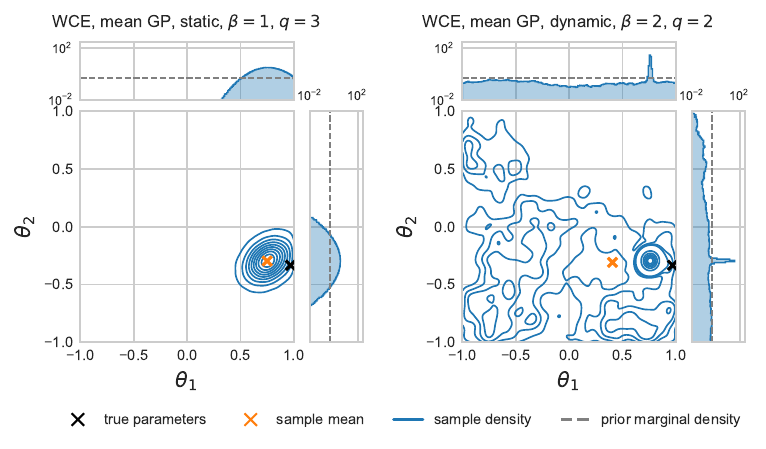}
    \caption{Joint posterior density contours of $\btheta$ with marginal histograms (log-frequency) of WCE posterior samples with the static weighting (left), and the dynamic weighting (right).}    
    \label{fig:wce_static_dynamic}
\end{figure}

The weighted cross-entropy can address bias-variance trade-off and improve inference via an observation dependent static scaling with an adjusted translation invariant weight kernel such as the IMQ weight in \eqref{eq:w_squared} for a reasonable choice of tuning parameter. When considering a parameter dependent scaling, the shape of the posterior may be significantly altered. The IMQ weight in dynamic weighting shrinks both for a wrong parameter and for a wrong surrogate, the residual-based weight cannot distinguish the two. Driven by $q^2$, the dynamic WCE posterior is a mixture of the prior and the CE posterior without a regime of moderate inflation around the mode. This is best understood via Figure~\ref{fig:wce_static_dynamic}, the posterior maintains the concentrated mode of the mean GP surrogate CE posterior and quickly recovers the prior shape, here uniform, away from the mode. If the true parameter is not within the peak of the mode, the dynamic WCE posterior assigns it the prior density, hence the values in Table~\ref{tab:results_NLL_90}. Furthermore, since $w^2(\by,\btheta)\,\mathrm{QF}(\by,\btheta)\leq q^2$ for dynamic weighting, the corresponding WCE posterior with IMQ weight can differ from the prior by at most $\mathrm{KL}\left[p_0(\btheta) \,\|\, \widetilde{p}_\beta(\btheta \,|\, \by)\right]\leq\frac{\beta q^2}{2}$ ($q^2\to0$ always maintains the prior) for a single observation, suggesting an alternative specification of dynamic weights. Additionally, the shape of the dynamic WCE posterior was observed to be challenging for the MALA sampler.

The learning rate $\beta=2$ was motivated by counteracting overinflation of the marginal GP surrogate and by considering consistency ($\beta w^2\approx1$). For $n=5$ training points, the opposite applies since the CE posterior with the marginal GP surrogate is already overconfident. The additive covariance underestimates the surrogate error and the increased learning rate $\beta=2$ lowers coverage further ($79\%$ coverage instead of $90\%$ for the static weighting in Table~\ref{tab:results_NLL_90}). In NLPD, the static weighting differs little between $\beta=1$ and $\beta=2$. For an increasing number of training points $n$ and vanishing approximation error, we expect this to change.   

The tuning heuristic in Lemma~\ref{lem:tuning} connects to this. Already limited in applicability, it assumes the marginal GP well-specified which is not given for $n=5$ training points. Additionally, the Jensen gap can be considered to affect accuracy for $d_\by=6$. Yet, it is a first approach to tuning when numerical tuning is not available via expanding on \citet{reimann2026kalman} and can likely be adjusted to consider static weights. In the simulation experiment, it is overly conservative for the static weighting as indicated by the coverage in Table~\ref{tab:results_NLL_90} and by comparing the default tuning to Table~\ref{tab:grid_q}.

The presented approach made direct use of available tools for fast posterior evaluation via GP surrogates  in \citet{bai2024gaussian} and weighted cross-entropy posteriors in \citet{duran2024outlier} and explored to combine both. Section~\ref{sec:the_prop} analytically investigated the trade-off between the limits of robustness and consistency and Section~\ref{sec:sampling_methodology} considered implementation via standard methods maintaining tractability. Each aspect, modelling and design, theoretical analysis, and practice of inference, poses several follow-up question as discussed, e.g. weight design in dynamic WCE. We do not want to claim that generalized Bayesian inference is necessarily the right framework to address overfitting to the approximation bias. Yet, GBI enabled the results presented for an investigation on mitigating surrogate approximation error in the inference step in Bayesian inverse problems.

\newpage

\begin{funding}
Hans Reimann was supported by the Graduate Academy and the Heidelberg Graduate School of Mathematical and Computational Methods for the Sciences (HGS MathComp) at Heidelberg University. Aretha L Teckentrup was supported by EPSRC grants EP/X01259X/1 and EP/Y028783/1.
\end{funding}

\begin{disclosure}
    \vspace{-0.25cm}
    \paragraph{Programming.}
    The forward model, surrogates, samplers and evaluation metrics are implemented in Python~3.12.13 (NumPy~2.4.5, SciPy~1.17.1, scikit-learn~1.8.0, and Numba~0.65.1). The elliptic PDE is solved via finite elements in DOLFINx~0.9.0 (FEniCSx) with a direct LU solve through PETSc~3.24.3. The Python environment runs in a Docker container (Debian~13, Docker~29.4.3, WSL2). Simulation experiments were orchestrated via R~4.5.1. Monte Carlo repetitions were run in parallel with the base package \texttt{parallel} (\texttt{parLapply}, \texttt{parLapplyLB}) and executed the Pythons scripts via \texttt{callr}~3.7.6. The code is available on request since the implementation is standard and the changes via the proposed method are elementary.
    \vspace{-0.25cm}
    \paragraph{Generative AI.} Claude Opus 5.5 by Anthropic was used for the following tasks between September and November 2026 to revise a previous draft.
    \begin{itemize}
    \item \emph{Review.} The previous draft was reviewed for theoretical inconsistencies, ambiguities and oversights, and for gaps between the text and the code.
    \item \emph{Code.} The previous code base was reworked to improve suit orchestration and parallelisation in \texttt{R}, the surrogate functions were adjusted to match changes in the paper, and tuning, evaluation metrics, function tests and data visualization were implemented.
    \item \emph{Theoretical results.} Theorem 3.4 was adjusted to drop the observation dependence of the constant and Corollary 3.7 was restated in terms of the largest eigenvalue of the GP covariance.
    \item \emph{Language editing.} Spelling, grammar, notation and ambiguous wording were flagged and revisions suggested.
    \end{itemize}
    The authors checked all AI-assisted output: proofs were reworked by the authors and code was validated with test cases. Changes to the manuscript were made only by the authors. 
\end{disclosure}


\bibliographystyle{ba2}
\bibliography{refrences}

\end{document}

%% file: results/grid_q.tex
\begin{table}[tb]
\centering
\begin{tabular}{lcccc}
\toprule
Surrogate & static, $\beta = 1$ & dynamic, $\beta = 1$ & static, $\beta = 2$ & dynamic, $\beta = 2$\\
\midrule
mean GP & 3 & 0.5 & 2 & 2\\
marginal GP & 5 & 2 & 5 & 2\\
\bottomrule
\end{tabular}
\caption{\label{tab:grid_q}Grid-tuned scaling factor $\widetilde{q}$ per configuration ($q^2=\widetilde{q}^2\,d_\by$) via lowest average NLPD over $M_{\mathrm{train}}=50$ repetitions.}
\end{table}

%% file: results/table_NLL_90.tex
\begin{table}[tb]
  \centering
  \setlength{\tabcolsep}{4pt} 
  \begin{tabular}{@{}llcccc@{}}
    \toprule
    & & \multicolumn{2}{c}{Negative log-posterior} & \multicolumn{2}{c}{$90\%$ CI coverage (\%)} \\
    \cmidrule(lr){3-4} \cmidrule(lr){5-6}
    Surrogate & Posterior & default $q$ & tuned $q$ & default $q$ & tuned $q$ \\
    \midrule
    mean GP & CE & \multicolumn{2}{c}{3.02 (1.54)} & \multicolumn{2}{c}{4} \\
     & static, $\beta = 1$ & -0.453 (0.693) & \textbf{-1.93} (1.15) & 97 & 86 \\
     & dynamic, $\beta = 1$ & 1.39 (0.0981) & 1.39 (0.101) & 92 & \textbf{92} \\
     & static, $\beta = 2$ & -0.966 (0.689) & \textbf{-1.95} (1.08) & 97 & \textbf{88} \\
     & dynamic, $\beta = 2$ & 1.39 (0.0964) & 1.88 (1.35) & 92 & 48 \\
    \addlinespace
    marginal GP & CE & \multicolumn{2}{c}{-0.496 (2.57)} & \multicolumn{2}{c}{38} \\
     & static, $\beta = 1$ & -0.0996 (0.741) & \textbf{-2.07} (0.818) & 97 & \textbf{90} \\
     & dynamic, $\beta = 1$ & 0.415 (0.836) & -1.52 (1.82) & 96 & 57 \\
     & static, $\beta = 2$ & -0.620 (0.788) & \textbf{-2.0} (1.13) & 98 & 79 \\
     & dynamic, $\beta = 2$ & 0.136 (1.19) & -0.0370 (2.52) & 96 & 31 \\
    \addlinespace
    exact (FEM) & -- & \multicolumn{2}{c}{-8.56 (0.813)} & \multicolumn{2}{c}{91} \\
    prior & -- & \multicolumn{2}{c}{1.39 (0.0718)} & \multicolumn{2}{c}{92} \\
    \bottomrule
  \end{tabular}
  \caption{\label{tab:results_NLL_90}Negative log-posterior density at the true parameter (mean and empirical standard deviation) and coverage of the central $90\%$ marginal credible intervals over $100$ repetitions for default and grid-tuned $q$.}
\end{table}